\documentclass[10pt]{elsart}
\usepackage{amssymb}

\def\cal#1{\mathcal{#1}}
\def\a{\alpha}
\def\b{\beta}
\def\t{\tau}

\def\vectorfields#1{{\cal X}(#1)}

\def\G{\Gamma}

\def\lie#1{{\cal L}_{#1}}
\def\fpd#1#2{\frac{\partial #1}{\partial #2}}
\def\R{{\rm I\kern-.20em R}}
\def\del{\nabla}
\def\ovl#1{\overline{#1}}

\def\im{\mbox{Im }}

\def\half{\mbox{$\frac{1}{2}$}}

\def\beq{\begin{equation}}
\def\eeq{\end{equation}}
\def\bea{\begin{eqnarray*}}
\def\eea{\end{eqnarray*}}
\def\qed{\hskip0.1em\hfill\null\ \null\nobreak\hfill
\kern3pt\lower1.8pt\vbox{\hrule\hbox
{\vrule\kern1pt\vbox{\kern1.7pt \hbox{$\scriptstyle
QED$}\kern0.2pt}\kern1pt\vrule}\hrule}}

\def\oneforms#1{{\cal X}^*(#1)}

\def\ro{\rho}

\def\een{\mbox{id}}

\begin{document}
\begin{frontmatter}
\title{A connection theoretic approach to sub-Riemannian geometry}
\author{B. Langerock}
\address{
Department of Mathematical Physics and Astronomy, Ghent
University, Krijgslaan 281 S9,B-9000 Gent Belgium}

\begin{abstract}
We use the notion of generalized connection over a bundle map in
order to present an alternative approach to sub-Riemannian
geometry. Known concepts, such as normal and abnormal extremals,
will be studied in terms of this new formalism. In particular,
some necessary and sufficient conditions for the existence of
abnormal extremals will be derived. We also treat the problem of
characterizing those curves that verify both the nonholonomic
equations and the so-called vakonomic equations for a `free'
particle submitted to some kinematical constraints.
\end{abstract}

\end{frontmatter}
\section{Introduction} \label{intro}
A sub-Riemannian structure on a manifold is a generalization of a
Riemannian structure in that a metric is only defined on a proper
vector sub-bundle of the tangent bundle to the manifold (i.e.\ on
a regular distribution), rather than on the whole tangent bundle.
As a result, in sub-Riemannian geometry a notion of length can
only be assigned to a certain privileged set of curves, namely
curves that are tangent to the given regular distribution on which
the metric is defined. The problem then arises to find those
curves that minimize length, among all curves connecting two given
points. The characterization of these length minimizing curves is
one of the main research topics in sub-Riemannian geometry, which
has also interesting links to control theory and to vakonomic
dynamics (for the latter, see for instance J. Cort\'es, \etal{}
\cite{manolo}).

The connection with control theory will be touched upon in Section
\ref{gendef} where, without entering into the details, we will
present a formulation of the Maximum Principle, following the work
of R.S. Strichartz \cite{stri,stri2} and H.J. Sussmann
\cite{sussmann4}. This will lead, among others, to the definition
of \emph{normal} and \emph{abnormal extremals}. The connection
with vakonomic dynamics will be explored in Section
\ref{sec:VakonomicDynamicsAndNonholonomicMechanics}.

The main goal of this paper is to give an application to
sub-Riemannian geometry of the theory of generalized connections
over a bundle map, developed in a previous paper in collaboration
with F. Cantrijn \cite{onszelf}. In Section
\ref{sec:connectionsonasub-rstructure} we consider some aspects of
this theory in the framework of sub-Riemannian geometry. Then,
normal extremals will appear as ``geodesics'' and abnormal
extremals as ``base curves of parallel transported sections" with
respect to a suitable generalized connection associated to the
sub-Riemannian structure. Apart from shedding some new light on
certain elements of sub-Riemannian geometry, this formulation also
allows us to prove some known results in an elegant way.

The main subtlety in studying length minimizing curves of a
sub-Riemannian structure lies in the existence of ``abnormal
minimizers'', i.e. length minimizing abnormal extremals. R.
Montgomery was the first to construct an explicit example of such
abnormal curves (see \cite{Mont}). Since then, many other examples
were found, for instance by W. Liu and H.J. Sussmann in
\cite{sussmann3}. We will deal with this topic in Section
\ref{sec:AbnormalExtremals}, where necessary and sufficient
conditions for the existence of abnormal extremals are given.

In this paper, we only consider real, Hausdorff, second countable
smooth manifolds, and by smooth we will always mean
$C^{\scriptscriptstyle\infty}$. The set of (real valued) smooth
functions on a manifold $M$ will be denoted by $\cal{F}(M)$, the
set of smooth vector fields by $\vectorfields{M}$ and the set of
smooth one-forms by $\oneforms{M}$. Let $V$ be a real vector
space, and $W$ a subspace, then the annihilator space of $W$ is
given by
\[W^0 = \{\b \in V^* \; | \; \langle \b, w \rangle =0 \; \forall w
\in W\}.\] If $E$ is a vector bundle over a manifold $M$ and $F$
any vector sub-bundle, then the annihilator bundle $F^0$ of $F$ is
the sub-bundle of the dual bundle $E^*$ of $E$ over $M$ whose
fibre over a point $x \in M$ is the annihilator space of the
subspace $F_x$ of $E_x$. The set of smooth (local) sections of an
arbitrary bundle $E$ over $M$ is denoted by $\G(E)$. In this
paper, the domain of a curve will usually be taken to be a closed
(compact) interval in $\R$. Whenever we say that such a curve,
defined on an interval $[a,b]$, is an integral curve of a vector
field, we simply mean that it is the restriction of a maximal
integral curve defined on an open interval containing $[a,b]$.

\section{General Definitions}\label{gendef}
In this section, we first give a brief review of some natural
objects associated to a sub-Riemannian structure and we recall the
necessary conditions, derived from the Maximum Principle, for a
curve to be length minimizing. Next, we discuss some general
aspects of the theory of connections over a bundle map.

\subsection{sub-Riemannian structures: preliminary definitions}
\label{sbs:subRiemannianStructuresPreliminaryDefinitions} Suppose
that $M$ is a smooth manifold of dimension $n$, equipped with a
regular distribution $Q \subset TM$ (i.e. $Q$ is a smooth
distribution of constant rank, say of rank $k$). In view of the
regularity, $Q$ can alternatively be regarded as a vector
sub-bundle of $TM$ over $M$. The natural injection $i: Q
\hookrightarrow TM$ is then a linear bundle mapping fibred over
the identity. A regular distribution is also completely
characterized by its annihilator, i.e. giving $Q$ is equivalent to
specifying the sub-bundle $Q^0$ of the cotangent bundle $T^*M$
whose fibre over $x \in M$ consists of all co-vectors at $x$ which
annihilate all vectors in the subspace $Q_x$ of $T_xM$.

A smooth Riemannian bundle metric $h$ on $Q$ is a smooth section
of the tensor bundle $Q^* \otimes Q^* \rightarrow M$ such that it
is symmetric and positive definite, i.e. for all $X_x, Y_x \in
Q_x$ one has
\begin{eqnarray*}
h(x)(X_x, Y_x) &=& h(x)(Y_x,X_x),\\
h(x)(X_x,X_x) &\geq& 0, \mbox{ and the equality holds iff } X_x=0.
\end{eqnarray*}
With a Riemannian bundle metric one can associate a smooth linear
bundle isomorphism $\flat_h: Q \rightarrow Q^*,\; X_x \mapsto
h(x)(X_x,.)$, fibred over the identity on $M$, with inverse
denoted by $\sharp_h:=\flat_h^{-1}: Q^* \rightarrow Q$.
\begin{defn}
A sub-Riemannian structure $(M,Q,h)$ is a triple where $M$ is a
smooth manifold, $Q$ a smooth regular distribution on $M$, and $h$
a Riemannian bundle metric on $Q$.
\end{defn}
Although it is not explicitly mentioned in the definition, it will
always be tacitly assumed, as it is customary in sub-Riemannian
geometry, that $Q$ is a non-integrable distribution and,
therefore, does not induce a foliation of $M$. A manifold $M$
equipped with a sub-Riemannian structure, will be called a \emph{
sub-Riemannian manifold}. With a sub-Riemannian structure
$(M,Q,h)$ one can associate a smooth mapping $g: T^*M \rightarrow
TM$ defined by
\[g(\a_x) = i \left( \sharp_h \left(i^*(\a_x) \right) \right) \in
TM,
\]
where $i^*: T^*M \rightarrow Q^*$ is the adjoint mapping of $i$,
i.e. for any $\a_x \in T^*_xM$, $i^*(\a_x)$ is determined by
$\langle i^*(\a_x),X_x\rangle = \langle \a_x,i(X_x)\rangle$, for
all $X_x \in Q_x$. Clearly, $g$ is a linear bundle mapping whose
image set is precisely the sub-bundle $Q$ of $TM$ and whose kernel
is the annihilator $Q^0$ of $Q$. To simplify notations we shall
often identify an arbitrary vector in $Q$ with its image in $TM$
under $i$ and smooth sections of $Q$ (i.e. elements of $\G(Q)$)
will often be regarded as vector fields on $M$.

With $g$ we can further associate a section $\ovl g$ of $TM
\otimes TM \rightarrow M$ according to
\[\ovl g(x)(\a_x, \b_x) =
\langle g(\a_x), \b_x \rangle
\]
for all $x \in M$ and $\a_x, \b_x \in T^*_xM$. From \bea \ovl
g(x)(\a_x,\b_x) :=
\langle g(\a_x) , \b_x \rangle &=& \langle \sharp_h (i^* \a_x), i^*(\b_x) \rangle\\
& = & h(x)(\sharp_h(i^* \a_x) , \sharp_h(i^*\b_x))\\
& = & h(x)(g(\a_x),g(\b_x)), \eea we conclude that $\ovl g$ is
symmetric.

Let $G$ be a Riemannian metric on $M$. It is easily seen that,
given a regular distribution $Q$ on $M$, we can associate with the
metric $G$ a sub-Riemannian structure $(M,Q,h_G)$ where $h_G$ is
the restriction of $G$ to the sub-bundle $Q$, i.e.
$h_G(x)(X_x,Y_x) := G(x)(X_x, Y_x)$ for any $x \in M$ and $X_x,Y_x
\in Q_x$. Given a sub-Riemannian structure $(M,Q,h)$ and a
Riemannian metric $G$ on $M$, we say that the Riemannian metric
restricts to $h$ if $h_G = h$. Now, every sub-Riemannian structure
can be seen as being determined (in a non-unique way) by the
restriction of a Riemannian metric. Indeed, let $h$ be a
Riemannian bundle metric on a vector sub-bundle $Q$ of $TM$, and
let $\{U_\a\}$ be an open covering of $M$ such that, on each
$U_{\a}$, there exists an orthogonal basis $\{X_1, \ldots , X_k\}$
of local sections of $Q$ with respect to $h$. Extend this to a
basis of vector fields $\{X_1, \ldots , X_n\}$ on $U_{\a}$ and
define a Riemannian metric on $U_{\a}$ by
\[G_{\a}(x)(X_x,Y_x) = \sum_{i,j=1}^k a^ib^jh(x)(X_i(x),X_j(x)) +
\sum_{i=k+1}^n a^ib^i,\]where $X_x=a^iX_i(x)$ and $Y_x=b^iX_i(x)$,
with $a^i,b^i \in \R$.  One can then glue these metrics together,
using a partition of unity subordinate to the given covering
$\{U_{\a}\}$. This procedure, which is similar to the one adopted
for constructing a Riemannian metric on an arbitrary smooth
manifold (see for instance \cite{Brickell}, Proposition 9.4.1),
produces a Riemannian metric on $M$ which, by construction,
restricts to $h$.

In the sequel we will repeatedly make use of a Riemannian metric
$G$ which restricts to a given sub-Riemannian metric $h$. In that
connection we now introduce some further notations and prove some
useful relations associated to $G$ and $h$. The natural bundle
isomorphism between $TM$ and $T^*M$ induced by $G$ will be denoted
by $\sharp_G$, with inverse $\flat_G = \sharp_G^{-1}$. Let $x \in
M$ and let $X_x,Y_x \in Q_x$, then one has:
\[\langle i^*\flat_G(i(X_x)), Y_x \rangle =\langle
\flat_G(i(X_x)), i(Y_x) \rangle = \langle \flat_h(X_x), Y_x
\rangle\;,\] which implies that $\flat_h = i^*\circ\flat_G\circ
i$. Inserting this into $g \circ \flat_G \circ i$ and taking into
account the definition of $g$, we conclude that \[g \circ \flat_G
\circ i = i \mbox{ or } g \circ \flat_G|_Q = \een_Q\;,\]where
$\een_Q$ is the identity mapping on $Q$. The orthogonal
projections of $TM$ onto $Q$ and onto its $G$-orthogonal
complement $Q^\bot$ will be denoted by $\pi$ and $\pi^\bot$,
respectively. Now, $T^*M$ can be written as the direct sum of
$(Q^{\bot})^0$ and $Q^0$ and the corresponding projections will be
denoted by $\t$ and $\t^\bot$, respectively. It is easily proven
that $(Q^{\bot})^0 \cong \flat_G(Q)$ and that \[\t^\bot = \flat_G
\circ \pi^\bot \circ \sharp_G\;, \quad \t = \flat_G \circ \pi
\circ \sharp_G\;.\] Using the fact that $\left.g \circ \flat_G
\right|_Q = \een_Q$ and $\ker g= Q^0$, we also have: $g= g \circ
\t = \pi \circ \sharp_G$.

To any regular distribution $Q$ on $M$ one can associate a natural
tensor field acting on $Q^0 \otimes Q \otimes Q$. Indeed, let
$\eta \in \G(Q^0), X,Y \in\G(Q)$ and let $[X,Y]$ denote the Lie
bracket of $X$ and $Y$, regarded as vector field on $M$. Then it
is easily proven that the expression $\langle \eta, [X,Y]\rangle$
is $\cal{F}(M)$-linear in all three arguments and, therefore,
determines a tensorial object. Now, $Q$ is involutive if and only
if this tensor is identically zero. Next, assume that $Y \in
\vectorfields{M}$, with $\eta$ and $X$ as before, then $\langle
\eta, [X,Y]\rangle$ is still $\cal{F}(M)$-linear in $\eta$ and $X$
(but not in $Y$). This justifies the following notation, which
will be used later on in our discussion of the Maximum Principle:
for any $x \in M$, $\eta_x \in Q^0_x, X_x \in Q_x$ and arbitrary
$Y \in \vectorfields{M}$, put \begin{equation}\label{bracket}
\langle \eta_x,[X_x,Y]\rangle := \langle \eta,
[X,Y]\rangle(x)\;,\end{equation} where $\eta$ (resp. $X$) may be
any section of $Q^0$ (resp. $Q$) such that $\eta(x) = \eta_x$
(resp. $X(x)=X_x$).
\subsection{Necessary conditions for length minimizing curves}
\label{sbs:NormalAndAbnormalExtremalsLocalEquations} For the
further discussion in this paper it is important that we give a
precise description of the class of curves we will be dealing
with. First of all, by a \emph{curve} in an arbitrary manifold $P$
we shall always mean a smooth mapping (in the $C^{\infty}$ sense)
$c: I \rightarrow P$, with $I \subset \R$ a closed interval, and
such that $c$ admits a smooth extension to an open interval
containing $I$. A mapping $c:[a,b] \rightarrow P$ will be called a
\emph{piecewise curve} in $P$ if there exists a finite subdivision
$a_1:=a < a_2 < \ldots < a_k < a_{k+1}: =b$ such that the
following conditions are fulfilled: \begin{enumerate} \item $c$ is
left continuous at each point $a_i$ for $i =2, \ldots, k+1$, i.e.
$\lim_{t \to a_i^{-}}c(t)$ exists and equals $c(a_i)$;
\item $\lim_{t \to a_i^{+}}c(t)$ is defined for all $i = 1,
\ldots, k$ and $\lim_{t \to a_1^{+}}c(t)=c(a_1)$ (i.e.\ $c$ is
right continuous at $a_1=a$);
\item for each $i = 1,\ldots,k$, the mapping $c^i: [a_i,a_{i+1}]
\rightarrow P$, defined by $c^i(t) = c(t)$ for $t \in
]a_i,a_{i+1}]$ and $c^i(a_i) = \lim_{t \to a_i^{+}}c(t)$, is
smooth (i.e. is a curve in $P$).
\end{enumerate} A piecewise curve which is continuous everywhere,
will simply be called a \emph{continuous piecewise curve} (and
corresponds to what is often called in the literature, a piecewise
smooth curve.)

In the sequel, whenever we are dealing with a (continuous)
piecewise curve $c: [a,b] \rightarrow P$, the notation $c^i$ will
always refer to the curve defined on the $i^{th}$ subinterval of
$[a,b]$, bounded by points where $c$ fails to be smooth.

Consider now a sub-Riemannian structure $(M,Q,h)$, with associated
bundle map $g: T^*M \rightarrow TM$. A curve (resp.\ piecewise
curve) $c: [a,b] \rightarrow M$ is said to be \emph{tangent to
$Q$} if $\dot{c}(t) \in Q_{c(t)}$ for all $t \in [a,b]$ (resp. for
all $t$ where the derivative exists). Next, let $\alpha: I
\rightarrow T^*M$ be a curve in $T^*M$ and put $c = \pi_M \circ
\alpha$, with $\pi_M: T^*M \rightarrow M$ the natural cotangent
bundle projection. Then, we say that $\alpha$ is
\emph{$g$-admissible} if
\[ g(\alpha(t)) = \dot{c}(t),\quad \hbox{for all}\; t \in I\;.\]
The projected curve $c$ will be called the \emph{base curve} of
$\alpha$. If $\alpha: I = [a,b] \rightarrow T^*M$ is a piecewise
curve, then $\alpha$ will be called $g$-admissible if its
projection $c = \pi_M \circ \alpha$ is a continuous piecewise
curve such that, in addition, $g(\alpha^i(t)) = \dot{c}^i(t)$ for
$t \in [a_i,a_{i+1}]$ (where we have used the notational
conventions introduced above). We now prove the following result
which will be of use later on.

\begin{lem}\label{g-admissible} Given a sub-Riemannian structure
$(M,Q,h)$ and any curve (resp.\ continuous piecewise curve) $c$ in
$M$, tangent to $Q$. Then, there always exists a $g$-admissible
curve (resp.\ piecewise curve) in $T^*M$ which projects onto $c$.
\end{lem}
\begin{pf} Take a Riemannian metric $G$ which restricts to $h$ on $Q$.
If $c: [a,b] \rightarrow M$ is a curve tangent to $Q$, one can
simply put $\alpha(t) = \flat_G(\dot{c}(t)$ for all $t \in [a,b]$.
Clearly, $\alpha$ then defines a $g$-admissible curve in $T^*M$
with base curve $c$.

Next, assume $c: [a,b] \rightarrow M$ is a continuous piecewise
curve, tangent to $Q$. We can then define a piecewise curve
$\alpha$ in $T^*M$ as follows: put $\alpha(t) =
\flat_G(\dot{c}(t))$ for all $t$ where $\dot{c}(t)$ is defined
and, using the notational conventions introduced above,
$\alpha(a_{i+1}) = \flat_G(\dot{c}^i(a_{i+1}))$ for
$i=1,\ldots,k$. It is easy to check that the mapping $\alpha:
[a,b] \rightarrow T^*M$ thus constructed, is a $g$-admissible
piecewise curve, projecting onto $c$. \qed \end{pf}

We will now introduce the notion of length of curves, and of
continuous piecewise curves, tangent to $Q$.
\begin{defn}
Given a sub-Riemannian structure $(M,Q,h)$, then the length of a
curve $c: [a,b] \rightarrow M$, tangent to $Q$, is given by \[L(c)
:= \int_a^b \sqrt{h(c(t))(\dot c(t), \dot c(t))} dt.\]
\end{defn}
Given any $g$-admissible curve $\a$ in $T^*M$ with base curve $c$,
and a Riemannian metric $G$ which restricts to $h$, then the
length of $c$ still equals
\[L(c)= \int^b_a \sqrt{\ovl g(c(t))(\a(t) ,\a(t))} dt = \int^b_a
\sqrt{G(c(t))(\dot c(t), \dot c(t))} dt .\] In particular, the
value of these integrals do not depend on the specific choice of
$\a$, resp.\ $G$.

The above notion of length can be easily extended to the class of
continuous piecewise curves $c$, tangent to $Q$, by putting $L(c)
= \sum_{i=1}^{k}L(c^i)$.

For the following discussion, which is partially inspired on
Sussmann's approach to a coordinate free version of the Pontryagin
Maximum Principle \cite{sussmann4}, we make two additional
assumptions. First, we assume that $M$ is pathwise connected, and
secondly, we take the distribution $Q$ of the given sub-Riemannian
structure $(M,Q,h)$ to be bracket generating, i.e.\ if $L(Q)$
denotes the Lie algebra generated by sections of $Q$, regarded as
vector fields on $M$, then we assume that at each point $x \in M$,
$T_xM = \{X(x)\;|\; \hbox{for all}\; X \in L(Q)\}$. Both
assumptions imply in particular that any two points of $M$ can be
joined by a continuous piecewise curve tangent to $Q$, as follows
from a well-known theorem of Chow \cite{chow}. Therefore, under
these assumptions it makes sense to talk about the length
minimizing curve connecting two given points. More precisely,
given a continuous piecewise curve $c: [a,b] \rightarrow M$
tangent to $Q$, connecting two points $x_0$ and $x_1$ (i.e.\ $c(a)
= x_0$, $c(b)=x_1$), then $c$ is called \emph{length minimizing}
if $L(c)\le L(\tilde c)$ for any other continuous piecewise curve
$\tilde c :[a,b] \to M$ tangent to $Q$, with $\tilde c(a)= x_0$
and $\tilde c(b) = x_1$.

Note that, given a continuous piecewise curve $c$, connecting two
points $x_0$ and $x_1$, one can always determine a
parameterization of $c$ such that $c:[0,1] \rightarrow M$, with
$c(0)=x_0,c(1)=x_1$, and for which there exists a nonzero constant
$k$ such that $h(c(t))(\dot c(t), \dot c(t))=k$ for all $t$ where
$\dot{c}(t)$ is defined. Following Sussmann, we will call this a
\emph{parameterization by constant times arc-length}.

We now arrive at the following weak version of the Maximum
principle.
\begin{thm} \label{cfmp} Consider a sub-Riemannian structure
$(M,Q,h)$ with $M$ connected and $Q$ bracket generating. Let $c:
[0,1] \rightarrow M$ be a continuous piecewise curve which is
length minimizing, and parameterized by constant times arc-length.
Then, there exists a continuous piecewise curve $\psi : [0,1] \to
T^*M$ along $c$, i.e. $\pi_M(\psi(t))= c(t)$, which does not
intersect the zero section and
such that at least one of the following two conditions is satisfied:\\
(i) $\psi(t)$ is an integral curve of the Hamiltonian vector field
$X_H$ on $T^*M$, with Hamiltonian given by the smooth function
$H(\a_x) = \half \ovl g(x)(\a_x,\a_x)$ for $\a_x \in T_x^*M$,
which, in particular, implies that both $\psi$ and $c$ are smooth;\\
(ii) $\psi(t) \in Q^0$ for all $t \in I$, and for any piecewise
$g$-admissible curve $\a$ with base curve $c$, the following
equation holds:
\[\left.\frac{d}{dt}\right|_t \langle
\psi(t), X(c(t))\rangle = \langle \psi(t), [g(\a(t)), X]
\rangle,\]for all $X \in \vectorfields{M}$ and all $t \in [0,1]$
where $\dot c(t)$ is well defined.
\end{thm}

(Note that on the right-hand side of the equation in (ii) we have
used the notation introduced in (\ref{bracket})). For a derivation
of this weak version of the Maximum principle in terms of the more
general class of absolutely continuous curves, we refer to
\cite{sussmann4}. Inspired on the (local) analysis presented in
\cite{Pont} (p. 79), the proof of the above theorem follows by
making some minor adjustments to the one given in
\cite{sussmann4}.

\begin{defn}\label{normabnorm}
A continuous piecewise curve $c$ tangent to $Q$ is called a {\rm
normal} (resp. {\rm abnormal}) {\rm extremal} if there exists a
continuous piecewise curve $\psi$ in $T^*M$ along $c$, which does
not intersect the zero section of $T^*M$, satisfying condition (i)
(resp. (ii)) of Theorem \ref{cfmp}.
\end{defn}

Note that, according to this definition, normal or abnormal
extremals do not have to be length minimizing and that $c$ can be
simultaneously a normal and an abnormal extremal.

\subsection{Connections over a bundle map: general setting}
\label{sbs:ConnectionsOverABundleMapGeneralSetting} Inspired by
some recent work of R.L. Fernandes on ``contravariant
connections'' in Poisson geometry (see \cite{ferdi}) and, more
generally, connections associated with Lie algebroids (see
\cite{ferdi2}), we have recently embarked on the study of a
general notion of connection, namely connections defined over a
vector bundle map. This concept covers, besides the standard
notions of linear and nonlinear connections, various
generalizations such as partial connections and
pseudo-connections, as well as the Lie algebroid connections
considered by Fernandes. For a detailed treatment we refer to a
forthcoming paper, written in collaboration with F. Cantrijn
\cite{onszelf}. After briefly sketching the main idea underlying
the notion of a generalized connection over a vector bundle map,
we shall apply this notion of connection to a sub-Riemannian
structure.

Let $M$ be a manifold and $\nu: N \rightarrow M$ a vector bundle
over $M$. Assume, in addition, that a linear bundle map $\ro: N
\rightarrow TM$ is given such that $\t_M \circ \ro = \nu$, where
$\t_M:TM \rightarrow M$ denotes the natural tangent bundle
projection. Note that we do not require $\rho$ to be of constant
rank. Hence, the image set $\im \rho$ need not be a vector
sub-bundle of $TM$, but rather determines a generalized
distribution as defined by P. Stefan and H.J. Sussmann (see e.g.
\cite{Lib}, Appendix 3). It follows that $\rho$ induces a mapping
of sections, $\G(N) \rightarrow \vectorfields{M}: s \mapsto \ro
\circ s$, also denoted by $\rho$. Next, let $\pi:E \rightarrow M$
be an arbitrary fibre bundle over $M$. We may then consider the
pull-back bundle $\tilde{\pi}_1:\pi^{\ast}N \rightarrow E$, which
is a vector bundle over $E$. Note that $\pi^{\ast}N$ may also be
regarded as a fibre bundle over $N$, with projection denoted by
$\tilde{\pi}_2: \pi^{\ast}N \rightarrow N$.

\begin{defn}A generalized connection on $E$ over
the bundle map $\rho$ (or, shortly, a {\rm $\rho$-connection on
$E$}) is then defined as a linear bundle map $h: \pi^{\ast}N
\rightarrow TE$ from $\tilde{\pi}_1$ to $\tau_E$, over the
identity on $E$, such that, in addition, the following diagram is
commutative

\begin{picture}(1,3)(-3,2.2)
\thicklines \put(2,2.5){$N$} \put(2.7,2.7){\vector(1,0){2.5}}
\put(5.5,2.5){$TM$} \put(5.7,4.3){\vector(0,-1){1.3}}
\put(5.5,4.6){$TE$} \put(2.7,4.8){\vector(1,0){2.5}}
\put(1.7,4.6){${\pi}^{\ast}N$} \put(2.2,4.3){\vector(0,-1){1.3}}
\put(3.8,2.1){$\rho$} \put(3.8,5.2){$h$} \put(6,3.5){$T{\pi}$}
\put(1.5,3.5){$\tilde{\pi}_2$}
\end{picture}

(where $T\pi$ denotes the tangent map of $\pi$).\end{defn} The
image set $\im h$ determines a generalized distribution on $E$
which projects onto $\im \rho$. It is important to note that $\im
h$ may have nonzero intersection with the bundle $VE$ of
$\pi$-vertical tangent vectors to $E$. The standard notion of a
connection on $E$ is recovered when putting $N=TM$, $\nu =
\tau_M$, and $\rho$ the identity map. In case $P$ is a principal
$G$-bundle over $M$, with right action $R: P \times G \rightarrow
P, (e,g) \mapsto R(e,g) = R_g(e) (= eg)$, a $\rho$-connection $h$
on $P$ will be called a \emph{principal $\rho$-connection} if, in
addition, it satisfies
\[
TR_g(h(e,n)) = h(eg,n),
\]
for all $g \in G$ and $(e,n) \in \pi^{\ast}N$. Slightly modifying
the construction described by Kobayashi and Nomizu \cite{koba},
given a principal $\rho$-connection on $P$, one can construct a
$\rho$-connection on any associated fibre bundle $E$.

Assume $E$ is a vector bundle and let $\{\phi_t\}$ denote the flow
of the canonical dilation vector field on $E$. A $\rho$-connection
$h$ on $E$ is then called \emph{a linear $\rho$-connection\/} if
\[ T \phi_t(h(e,n)) = h(\phi_t(e),n),\]for all $(e,n) \in
\pi^{\ast}N$. In \cite{onszelf} it is shown that such a linear
$\rho$-connection can be characterized by a mapping $\del: \G(N)
\times \G(E) \rightarrow \G(E), (s,\sigma) \mapsto \del_s\sigma$
such that the following properties hold:
\begin{enumerate}
\item $\del$ is $\R$-linear in both arguments;
\item $\del$ is ${\mathcal F}(M)$-linear in $s$;
\item for any $f \in {\mathcal F}(M)$ and for all $s \in \G(N)$ and
$\sigma \in \G(E)$ one has: $\del_s (f \sigma) = f \del_s \sigma +
(\ro \circ s)(f)\sigma$.
\end{enumerate}
It immediately follows that $\del_s \sigma(m)$ only depends on the
value of $s$ at $m$, and therefore we may also write it as
$\del_{s(m)}\sigma$. Clearly, $\del$ plays the role of the
covariant derivative operator in the case of an ordinary linear
connection. Henceforth, we will also refer to the operator $\del$
as a linear $\rho$-connection. Let $k$ and $\ell$ denote the fibre
dimensions of $N$ and $E$, respectively, and let
$\{s^{\alpha}:\alpha = 1,\ldots,k\}$, resp. $\{\sigma^A: A=1,
\ldots,\ell\}$, be a local basis of sections of $\nu$, resp.
$\pi$, defined on a common open neighborhood $U \subset M$. We
then have $\del_{s^{\alpha}}\sigma^A = \G^{\a A}_{B}\sigma^B$, for
some functions $\G^{\a A}_{B}\in {\mathcal F}(U)$, called the
connection coefficients of the given $\rho$-connection. A
$\ro$-connection $\del$ can be extended to an operator, also
denoted by $\del$, acting on sections of any tensor product bundle
of $E$. This can be achieved by applying standard arguments, and
the details are left to the reader. We just like to mention here
that the action on ${\mathcal F}(M)$ and $\G(E^*)$ is determined
by the following relations: for $s \in \G(N)$, $f \in {\mathcal
F}(M)$, $\sigma \in \G(E)$ and $\zeta \in \G(E^*)$,
\[\del_sf := (\rho \circ s) (f)\;, \quad \del_s \langle \sigma, \zeta \rangle
=  \ro(s)\langle \sigma,\zeta \rangle = \langle \del_s \sigma ,
\zeta \rangle + \langle \sigma , \del_s \zeta \rangle.\]

In order to associate a notion of parallel transport to a linear
$\rho$-connection, we first need to introduce a special class of
curves in $N$. A curve $\tilde{c}: I=[a,b] \rightarrow N$ is
called \emph{$\ro$-admissible} if for all $t \in I$, one has
$\dot{c}(t) = (\ro \circ \tilde{c})(t)$, where $c = \nu \circ
\tilde{c}$ is the projected curve on $M$. Curves in $M$ that are
projections of $\ro$-admissible curves in $N$ are called
\emph{base curves}. (We will see that this terminology is in
agreement with the one introduced in the previous subsection.)
Note that, in principle, a base curve may reduce to a point.

As in standard connection theory, with any linear
$\rho$-connection $\del$ on a vector bundle $\pi: E \rightarrow
M$, and any $\ro$-admissible curve $\tilde{c}: [a,b] \rightarrow
N$, one can associate an operator $\del_{\tilde{c}}$, acting on
sections of $\pi$ defined along the base curve $c = \nu \circ
\tilde{c}$. The operator $\del_{\tilde{c}}$ is completely
determined by the following prescriptions. For arbitrary sections
$\sigma$ of $\pi$ along $c$ (i.e. curves $\sigma: [a,b]
\rightarrow E$, satisfying $\pi \circ \sigma = c$) and for
arbitrary $f \in {\cal F}([a,b])$:
\begin{enumerate}
\item $\del_{\tilde{c}}$ is $\R$ linear;
\item $\del_{\tilde{c}} f \sigma = \dot f \sigma  + f \del_{\tilde{c}} \sigma$;
\item $\del_{\tilde{c}} \sigma (t) = \del_{\tilde{c}(t)}
\overline{\sigma}$, for $\overline \sigma \in \G(\pi)$ such that
$\overline \sigma(c(t)) = \sigma(t)$ for all $t \in [a,b]$.
\end{enumerate}
\begin{defn}
A section $\sigma$ of $\pi$, defined along the base curve of a
$\ro$-admissible curve $\tilde{c}:[a,b] \rightarrow N$, will be
called {\rm parallel along $\tilde{c}$} if $\del_{\tilde{c}}\sigma
(t) = 0$ for all $t \in [a,b]$. \end{defn} Taking again
$\{s^{\alpha}\}$, resp. $\{\sigma^A\}$, to be a local basis of
sections of $\nu$, resp. $\pi$, and putting
$\sigma(t)=\sigma_A(t)\sigma^A(c(t))$ and
$\tilde{c}(t)=\tilde{c}_{\a}(t)s^{\a}(c(t))$, we find that
$\sigma$ is parallel along $\tilde{c}$ if
\[\del_{\tilde{c}}\sigma(t) = \left(\dot \sigma_A(t) + \G^{\a B}_
A(\tilde{c}(t)) \sigma_B(t) \tilde{c}_\a(t)\right) \sigma^A(c(t))
= 0,\] which gives a system of linear differential equations for
the components of $\sigma$. Again using standard arguments, one
can show that this leads to a notion of parallel transport on $E$
along $\ro$-admissible curves in $N$ (cf.\ \cite{onszelf} for more
details).

Suppose we are given two $\ro$-admissible curves $\tilde c^{i}:
[a_i,b_i] \rightarrow N$, $i=1,2$ with $\tilde{c}^1(b_1)$ and
$\tilde{c}^2(a_2)$ belonging to the same fibre of $\nu$, i.e.
$c^1(b_1) = c^2(a_2)$ (where $c^i$ is the base curve of
$\tilde{c}^i$). Given any point in $E_{c^{1}(a_1)}$ one can
construct a unique parallel section along $\tilde c^{1}$, starting
from that point. The endpoint of this curve (at $t= b_1$) lies in
the fibre $E_{c^{2}(a_2)}$ and, therefore, can be taken as the
initial point of a unique parallel curve along $\tilde c^{2}$.
This construction can now be easily extended to the class of
piecewise $\rho$-admissible curves defined below.

Recalling the definition of a piecewise curve, given in the
previous subsection, and using the notational conventions
introduced there, a \emph{piecewise $\ro$-admissible curve $\tilde
c$} is defined as a piecewise curve in $N$ such that: (i) for each
$i = 2,\ldots,k$, $\tilde{c}^i(a_i)$ and $\tilde{c}^{i+1}(a_i) (=
\lim_{t \to a_i^{+}}c(t))$ belong to the same fibre of $\nu$ or,
equivalently, the projection $c = \nu\circ \tilde c$ is a
continuous piecewise curve, (ii) $\ro(\tilde{c}^i(t)) =
\dot{c}^i(t)$ for all $i=1, \ldots, k$ and $t \in [a_i,a_{i+1}]$.
Extending the above construction in the case of two
$\rho$-admissible curves $\tilde{c}^1,\tilde{c}^2$, it is now
clear how to determine the notion of parallel transport along a
piecewise $\rho$-admissible curve.

The following class of linear $\rho$-connections will play an
important role in the further analysis.

\begin{defn} \label{defpartial}
A linear $\ro$-connection on a vector bundle $E$ is called {\rm
partial} if for any $\sigma \in \G(E)$ and $n \in \ker (\rho)$, we
have $\nabla_n\sigma = 0$.
\end{defn}
It is instructive to know that the condition for a connection to
be partial is equivalent to the property that no (nonzero)
vertical tangent vectors to $E$ exist that are also contained in
$\im h$, as stated in the following proposition. For the proof,
which is quite technical, we refer to \cite{onszelf}.
\begin{prop}
Let $\del$ be a linear $\ro$-connection. Then $\del$ is partial if
and only if $ \im h  \cap VE =\{0\}$.
\end{prop}
\section{Connections on a sub-Riemannian structure}
\label{sec:connectionsonasub-rstructure} Fix a sub-Riemannian
structure $(M,Q,h)$ and consider the associated bundle map $g :
T^*M \to TM$. In this section we will be interested in generalized
connections on $T^*M$ over $g$. Our main goal is the
characterization of normal and abnormal extremals of the
sub-Riemannian structure in terms of such generalized connections.
Let $U$ be the domain of a coordinate chart in $M$. We will always
denote coordinates on $U$ by $x^i$, $i=1, \ldots, n$. The
coordinates on the corresponding bundle chart of $T^*M$ are
denoted by $(x^i,p_i)$, $i=1, \ldots ,n$.

\begin{defn}
A $g$-connection on $(M,Q,h)$ is a linear generalized connection
on $T^*M$ over the bundle map $g: T^*M \rightarrow TM$.
\end{defn}
Comparing with the notations from the previous section, we see
that a $g$-connection on a sub-Riemannian manifold is a linear
$\ro$-connection with $N=E=T^*M$ and $\ro = g$. Note that with
these identifications, the definition of a $g$-admissible curve,
as given in the context of sub-Riemannian geometry, agrees with
the notion of a $\rho$-admissible curve.
\begin{defn}
A $g$-admissible curve $\a: I \to T^*M$ is said to be an {\rm
auto-parallel curve} with respect to a $g$-connection $\del$ if it
satisfies $\del_\a \a(t)=0$ for all $t \in I$. Its base curve $c =
\pi \circ \alpha$ is then called a {\rm geodesic} of $\nabla$.
\end{defn}

In coordinates, an auto-parallel curve $\a(t)= (x^i(t),p_i(t))$
satisfies the equations
\[\dot x^i(t)=g^{ij}(x(t))p_j(t)\;, \quad  \dot p_j(t) = -
\G^{ik}_j(x(t))p_i(t)p_k(t),\] where $g^{ij}$ and $\G^{ik}_j \in
{\cal F}(U)$ are the local components of the contravariant tensor
field $\ovl g$ associated to the sub-Riemannian structure (cf.
Subsection 2.1) and the connection coefficients of $\del$,
respectively. In fact, given a linear $g$-connection $\del$ one
can always define a smooth vector field $\G^\del$ on $T^*M$ whose
integral curves are auto-parallel curves with respect to $\nabla$.
In canonical coordinates, this vector field reads:
\[\G^\del(x,p) = g^{ij}(x)p_j \fpd{}{x^i} - \G^{ik}_j(x) p_ip_k
\fpd{}{p_j}.\] (A proof of this property follows by standard
arguments, and is left to the reader.) This implies, in
particular, that given any $\a_0 \in T^*M$, there exists an
auto-parallel curve $\a$ passing through $\a_0$. Note that it may
happen that two different auto-parallel curves correspond to the
same base curve (i.e.\ may project onto the same geodesic).

Now, we would like to find a $g$-connection on a sub-Riemannian
manifold whose geodesics are precisely the normal extremals.
Recalling the definition of a normal extremal (Definition
\ref{normabnorm}), it follows that we will have to look for a
$g$-connection $\del$ for which $\G^\del = X_H$, where $X_H$
denotes the Hamiltonian vector field corresponding to $H(\a_x)=
\half \ovl g(\a_x,\a_x) \in {\cal F}(T^*M)$. A first step in that
direction is the construction of a symmetric product associated
with a given $g$-connection, which fully characterizes the
geodesics of the $g$-connection under consideration.

Two linear $g$-connections $\del$ and $\ovl{\del}$ have the same
geodesics if and only the tensor field $D: \oneforms{M} \otimes
\oneforms{M} \rightarrow \oneforms{M}, (\a,\b) \mapsto \del_\a \b-
\ovl{\del}_\a\b$ is skew-symmetric, or equivalently $D(\a,\a)
\equiv 0$. In local coordinates, the components of $D$ are given
by $D^{ij}_k = \G^{ij}_k - \ovl \G^{ij}_k$, where $\G^{ij}_k$ and
$\ovl \G^{ij}_k$ are the connection coefficients of $\del$ and
$\ovl \del$, respectively. We immediately see that $D$ is
skew-symmetric iff $\G^\del = \G^{\ovl \del}$, proving the
previous statement. Define the \emph{symmetric product} of a
connection $\del$ as \[\langle \a : \b \rangle_\del := \del_\a \b
+ \del_\b \a\;, \quad \hbox{for}\; \a,\b \in \oneforms{M}\;.\]
(Observe that this is not a tensorial quantity, i.e.\ $\langle \a
: \b \rangle_\del$ is not ${\cal F}(M)$-linear in its arguments).
By replacing $\a$ by $\a +\b$ in $D(\a,\a)$ the following lemma is
easily proven.
\begin{lem}
The geodesics of a linear $g$-connection $\del$ are completely
determined by the symmetric product $\langle \a : \b \rangle_\del$
in the sense that, given two $g$-connections $\del$ and $\ovl
\del$, then both have the same geodesics if and only if $\langle
\a : \b \rangle_\del = \langle \a : \b \rangle_{\ovl \del}$, for
all $\a,\b \in \oneforms{M}$.
\end{lem}
In the following we shall construct a symmetric bracket of
$1$-forms, associated to a sub-Riemannian structure $(M,Q,h)$,
which coincides with the symmetric product of a $g$-connection
$\del$ on $T^*M$ iff $\G^\del = X_H$.

Before proceeding, we first recall that the Levi-Civita connection
$\del^G$ associated to an arbitrary Riemannian metric $G$ is
completely determined by the relation:
\bea2G(\del^G_XY, Z) &=& X(G(Y,Z)) + Y(G(X,Z)) -Z(G(X,Y)) \\
&&+ G([X,Y],Z) - G([X,Z],Y) - G(X,[Y,Z])\;,\eea for all $X,Y,Z \in
\vectorfields{M}$. This can still be rewritten as
\[2\flat_G(\del^G_XY) = \lie{X}\flat_G(Y) + \lie{Y}\flat_G(X) +
\flat_G([X,Y])-d(G(X,Y)),\] and the symmetric product of two
vector fields $X,Y$, defined by $\langle X : Y \rangle_{\del^G} =
\del^G_X Y + \del^G_Y X$, then satisfies
\[\flat_G (\langle X : Y \rangle_{\del^G}) = \lie{X}\flat_G(Y) +
\lie{Y}\flat_G(X) -d(G(X,Y)).\] The right-hand side of this
equation now inspires us to propose the following definition of a
symmetric bracket of $1$-forms on a sub-Riemannian manifold.
\begin{defn}
The symmetric bracket associated to a sub-Riemannian structure
$(M,Q,h)$ is the mapping $\{ \cdot , \cdot \} : \oneforms{M}
\times \oneforms{M} \to \oneforms{M}$ defined by:
\[\{\a,\b\} = \lie{g(\a)} \b + \lie{g(\b)}\a - d\left(\ovl g(\a,\b)\right).\]
\end{defn}
In the following proposition we list some properties of this
bracket, the first of which justifies the denomination ``symmetric
bracket''. The proofs of these properties are straightforward and
immediately follow from the above definition.
\begin{prop}
The symmetric bracket satisfies the following properties: for any
$\a,\b \in \oneforms{M}$
\begin{enumerate}
\item $\{\a,\b\} =\{\b,\a\}$;
\item the bracket is $\R$-bilinear;
\item $\{f\a,\b\} = g(\b)(f)\a + f \{\a,\b\}$, with $f\in {\cal F}(M)$,
\item $\{\a,\eta\} = \lie{g(\a)}\eta$, for any $\eta \in \G(Q^0)$, and
$\{\a,\eta\}=0$ if both $\a$ and $\eta$ belong to $\G(Q^0)$.
\end{enumerate} \end{prop}
The first three properties justify the following definition.
\begin{defn}
A $g$-connection $\del$ is said to be {\rm normal} if the
associated symmetric product equals the symmetric bracket, i.e. if
$\langle \a : \b \rangle_\del = \{\a,\b\}$ holds for all $\a,\b
\in \oneforms{M}$.
\end{defn}
The connection coefficients of a normal $g$-connection satisfy the
relations\[\G^{ij}_k+ \G^{ji}_k = \fpd{g^{ij}}{x^k}, \mbox{ for
all } i,j,k = 1, \ldots , n.\]

We are now going to introduce a special operator, determined by
the given distribution $Q$, which will play an important role
later on.

For that purpose, we first recall that, given a regular involutive
distribution $D$ on a manifold $M$, there exists a canonical
connection $\del^B$ on the bundle $D^0 \to M$ over the natural
injection $i: D \rightarrow TM$, sometimes called the `Bott
connection', defined by: $\del^B_X \eta = i_Xd\eta$, where $X \in
\G(D)$ and $\eta \in \G(D^0)$. Indeed, under the hypothesis that
$D$ is involutive, the image of $\del^B$ is again an element of
$\G(D^0)$. This connection was used by R. Bott in \cite{Bott} to
prove, among others, that certain Pontryagin classes of the bundle
$D^0\to M$ are identically zero. However, in the setting of a
sub-Riemannian structure $(M,Q,h)$, the distribution $Q$ is
assumed not to be involutive and, hence, the $1$-form $i_X d\eta$
in general will not belong to $\G(Q^0)$. Nevertheless, this
mapping naturally pops up in our approach to characterize normal
and abnormal extremals and, therefore, deserves some special
attention. More specifically, with any sub-Riemannian structure
$(M,Q,h)$ we can associate a mapping $\delta^B$ according to
\[\delta^B : \G(Q) \times \G(Q^0) \rightarrow
\oneforms{M}, (X,\eta) \mapsto \delta^B_X \eta = i_X d \eta.\]
(The superscript $B$ is kept to remind us of the fact that this
map reduces to the Bott connection in the case of involutive
distributions.)
\begin{defn}
Given a sub-Riemannian structure $(M,Q,h)$, a $g$-connection
$\del$ is said to be adapted to the bundle $Q$ (shortly {\rm
$Q$-adapted}) if $\del_\a \eta = \delta^B_{g(\a)} \eta$ for all
$\a \in \oneforms{M}$ and $\eta \in \G(Q^0)$.
\end{defn}
For the following theorem, recall the notation introduced in
Subsection \ref{sbs:subRiemannianStructuresPreliminaryDefinitions}
for the projection operators associated with a Riemannian metric
$G$ restricting to $h$, namely $\t: T^*M \to \flat_G(Q)$, $\t^\bot
: T^*M \to Q^0$.
\begin{thm} \label{thm1}
Let $\del$ be a $g$-connection, then the following statements are
equivalent:
\begin{enumerate}
\item \label{een} $\del$ is a normal $g$-connection;
\item \label{vier}
for all $\a \in \oneforms{M} : \del_\a \a = \half \{\a,\a\}$;
\item \label{drie}
$\langle \del_\a X, \b\rangle + \langle \del_\b X , \a \rangle =
\langle [g(\a),X] , \b \rangle + \langle [g(\b),X] , \a \rangle +
X(g(\a,\b))$ for all $\a,\b \in \oneforms{M}$ and $X \in
\vectorfields{M}$;
\item \label{twee} $\G^\del = X_H$ or,
equivalently, every geodesic of $\del$ is a normal extremal and
vice versa;
\item \label{vijf}
let $G$ be a Riemannian metric restricting to $h$ and let $\del^G$
be its Levi-Civita connection, then for all $\a \in \oneforms{M}$,
$\del$ satisfies:\[\del_\a \a = \del^G_{g(\a)}\t(\a) +
\delta^B_{g(\a)} \t^\bot(\a).\]
\end{enumerate} ({\rm Note that the right hand side of (\ref{drie})
agrees with the definition of the symmetrized covariant derivative
considered in \cite{stri}}).
\end{thm}
\begin{pf}
The equivalence of  {\it (\ref{een})} and {\it (\ref{vier})}
follows directly from the definition of a normal $g$-connection,
and the equivalence of {\it (\ref{een})} and {\it (\ref{drie})}
follows from $\langle \del_\a\b,X\rangle = g(\a)(\langle
\b,X\rangle) - \langle \b, \del_\a X \rangle$ after some tedious
but straightforward calculations.

{\it (\ref{vier})} $\Leftrightarrow$ {\it (\ref{twee})}. Choose an
arbitrary $\a_0 \in T^*M$. Let $U$ be a coordinate neighborhood of
$x_0 = \pi_M(\a_0)$ and put $\a_0= (x_0^i,p^0_j)$. Then,
$\del_\a\a = \half \{\a,\a\}$ implies, in particular, that the
connection coefficients $\G^{ij}_k$ of $\del$ on $U$ satisfy
\[\G^{ij}_k(x_0)p^0_ip^0_j = \frac{1}{2} \fpd{g^{ij}}{x^k}(x_0)p^0_i
p^0_j.\] The coordinate expression for the Hamiltonian vector
field $X_H$ at $\a_0$ equals: \[X_H (\a_0) = g^{ij}(x_0)p^0_j
\left.\fpd{}{x^i}\right|_{\a_0} - \frac{1}{2} \fpd{g^{ij}}{x^k}
p^0_i p^0_j \left.\fpd{}{p_k}\right|_{\a_0}.\] Recalling the
definition of $\G^\del$ it is easy to see that $\G^\del(\a_0) =
X_H(\a_0)$ for any $\a_0 \in T^*M$ if and only if $\del_\a \a =
\half \{\a,\a\}$ for each $\a \in \oneforms{M}$.

{\it (\ref{vier})} $\Leftrightarrow$ {\it (\ref{vijf})}. Let $G$
be a Riemannian metric restricting to $h$.  Recall the following
property of the Levi-Civita connection $\del^G$: \[\flat_G
(\langle X : Y \rangle_{\del^G}) = \lie{X}\flat_G(Y) +
\lie{Y}\flat_G(X) -d(G(X,Y)).\] Putting $X=Y=g(\a)$, this equation
becomes
\[\flat_G (\del^G_{g(\a)} g(\a)) =  \lie{g(\a)}\flat_G(g(\a)) - \half d (\ovl g(\a,\a)).\]
Using the identity $\flat_G(g(\a)) = \t(\a)$ derived in the
Subsection
\ref{sbs:subRiemannianStructuresPreliminaryDefinitions}, and
taking into account that $\del^G$ preserves the metric $G$, i.e.
$\del^G \circ \flat_G = \flat_G \circ \del^G$, we obtain
\bea\del^G_{g(\a)} \t(\a) &=& \lie{g(\a)} \t(\a) - \half d(\ovl
g(\a,\a)),
\\ &=& \half\{\a,\a\} - \lie{g(\a)} \t^\bot(\a).\eea
Since $\t^\bot(\a) \in \G(Q^0)$ and $g(\a)\in \G(Q)$, the last
term on the right-hand side reduces to $\delta^B_{g(\a)}
\t^\bot(\a)$, which completes the proof. \qed
\end{pf}
Theorem \ref{thm1} implies, in particular, that normal
$g$-connections exist. For instance, the mapping $\del$ defined by
$\del_\a \b = \del^G_{g(\a)} \t(\b) + \delta^B_{g(\a)}
\t^\bot(\b)$ is a linear $g$-connection and it is normal, in view
of the equivalence of {\it (\ref{een})} and {\it (\ref{vijf})}.
Moreover, for $\b \in \G(Q^0)$ we find that $\del_\a \b =
\delta^B_{g(\a)}\b$, i.e.\ the connection under consideration is
also $Q$-adapted. Summarizing, we have shown the following result.
\begin{prop}\label{existence} Given a sub-Riemannian structure
$(M,Q,h)$, one can always construct a normal and a $Q$-adapted
$g$-connection.\end{prop} Furthermore, the $g$-connection
constructed gives us a relation between a normal $g$-connection,
the Levi-Civita connection $\del^G$ of any Riemannian metric
restricting to $h$ and the operator $\delta^B$. This relation will
be very useful when we study the relation between vakonomic
dynamics and nonholonomic mechanics, (see Section
\ref{sec:VakonomicDynamicsAndNonholonomicMechanics}).

In the following theorem we shall characterize an abnormal
extremal in terms of a $Q$-adapted $g$-connection. According to
Definition \ref{normabnorm}, a continuous piecewise curve $c$
tangent to $Q$ is an abnormal extremal if there exists a
continuous piecewise section $\psi$ of $Q^0$ along $c$ such that
\begin{equation} \label{vgl1} \left.\frac{d}{dt}\right|_t \langle
\psi(t), X(c(t))\rangle = \langle \psi(t), [g(\a(t)), X] \rangle,
\end{equation} holds for an arbitrary chosen piecewise
$g$-admissible curve $\a$ projecting onto $c$, for any $X \in
\vectorfields{M}$ and for all $t$ where $\dot c(t)$ is defined. We
can now state the following interesting result.
\begin{thm}\label{thm2} Given a continuous piecewise curve $c: I \rightarrow M$,
tangent to $Q$. There exists a continuous piecewise section of
$Q^0$ along $c$ which is parallel with respect to a $Q$-adapted
$g$-connection if and only if $c$ is an abnormal extremal.
\end{thm}
\begin{pf}
Recall from Subsection
\ref{sbs:NormalAndAbnormalExtremalsLocalEquations} that a
continuous piecewise curve on $I=[a,b]$ is defined as a continuous
map which can be regarded as a concatenation of a finite number of
curves $c^i$ ($i= 1, \ldots,k)$, with domain, say $[a_i,a_{i+1}]
\subset I$ for $a_1=a <a_2 < \ldots < a_k < a_{k+1}=b$ and such
that $c^i(a_{i+1}) = c^{i+1}(a_{i+1})$.

Let $c$ be an abnormal extremal such that \ref{vgl1} holds. We
shall denote the curves associated to $\a$ and $\psi$ on the
subinterval $[a_i,a_{i+1}]$, by $\a^i$ and $\psi^i$, respectively.
Since $\psi$ is continuous, we have $\psi^i(a_{i+1}) =
\psi^{i+1}(a_{i+1})$. Then (\ref{vgl1}) can equivalently be
rewritten as:\[\left.\frac{d}{dt}\right|_t \langle \psi^i(t),
X(c^i(t))\rangle = \langle \psi^i(t), [g(\a^i(t)), X] \rangle,\;
\forall t \in [a_i,a_{i+1}] \; (i=1, \ldots ,k).\] Now, take a
$Q$-adapted $g$-connection $\del$ (which always exists in view of
Proposition \ref{existence}). By definition, $\del$ satisfies
$\del_\b \eta = \delta^B_{g(\b)} \eta$ for all $\b \in
\oneforms{M}$ and $\eta \in \G(Q^0)$. Now, assume $\del_\b \eta
=0$. This is clearly equivalent to the condition $\langle \del_\b
\eta, X \rangle = 0$ for any $X \in \vectorfields{M}$ which, in
view of the fact that $\del$ is $Q$-adapted, can be rewritten as
$\langle \lie{g(\b)} \eta , X \rangle =0$ or $g(\b)( \langle \eta,
X \rangle) = \langle \eta, [g(\b),X] \rangle$. Herewith, we have
proven that $\del_\b \eta = 0$ iff $g(\b)( \langle \eta, X
\rangle) = \langle \eta, [g(\b),X] \rangle$ for any $X \in
\vectorfields{M}$. This equivalence can be restated in the
following way. Given a $g$-admissible curve $\a^i$, with base
curve $c^i$ and $\psi^i$ a section of $Q^0$ along $c^i$, then
$\del_{\a^i} \psi^i(t) = 0$ if and only if
\[\frac{d}{dt} (\langle \psi^i(t), X(c^i(t))\rangle = \langle \psi^i(t),
[g(\a^i(t)),X]\rangle, \ \mbox{ for all } X \in
\vectorfields{M}.\] Now, $\del_{\a^i} \psi^i(t) = 0$ for all $t
\in [a_i,a_{i+1}]$ and $i=1, \ldots, k$, with $\psi^i(a_{i+1}) =
\psi^{i+1}(a_{i+1})$ implies, by definition, that the continuous
piecewise section $\psi$ of $Q^0$ is parallel with respect to the
$Q$-adapted $g$-connection $\del$ (see Section
\ref{sbs:ConnectionsOverABundleMapGeneralSetting}). This already
proves one half of the theorem.

The proof of the converse statement, namely that the existence of
a continuous piecewise section $\psi$ of $Q^0$ along $c$,
satisfying the appropriate conditions, implies that $c$ is an
abnormal extremal, simply follows by reversing the above
arguments.\qed
\end{pf}
To conclude this section we make some further remarks on normal
and $Q$-adapted $g$-connections. It is well known that the
Levi-Civita connection $\del^G$, associated with a Riemannian
metric $G$, is uniquely determined by the properties that it
preserves the metric, i.e. $\del^G G =0$, and that its torsion is
zero. We would like to consider now metric $g$-connections $\del$
on a sub-Riemannian manifold, i.e. $\del \ovl g=0$ (where $\ovl g$
is the symmetric contravariant $2$-tensor field defined in
Subsection
\ref{sbs:subRiemannianStructuresPreliminaryDefinitions}). From
above we know that normal extremals of a sub-Riemannian structure,
resp.\ abnormal extremals, can be characterized as geodesics of a
normal $g$-connection, resp.\ as parallel transported sections of
$Q^0$ for a $Q$-adapted $g$-connection (see Theorem \ref{thm1},
resp.\ Theorem \ref{thm2}). Therefore it is natural to look for
$g$-connections that are simultaneously normal and $Q$-adapted. It
has been shown above that such a $g$-connection always exists,
namely $\del_\a \b = \del^G_{g(\a)} \t(\b) + \delta^B_{g(\a)}
\t^0(\b)$, with $G$ any Riemannian metric restricting to $h$. We
will prove, however, that no metric $g$-connection can be found
that is also $Q$-adapted. First we prove an interesting result
relating the notion of partial $g$-connection (see Definition
\ref{defpartial}) with that of a $Q$-adapted normal
$g$-connection.
\begin{prop}
Let $\del$ be a normal $g$-connection. Then $\del$ is partial if
and only if $\del$ is $Q$-adapted.
\end{prop}
\begin{pf}
Let $\del$ be a normal $g$-connection, i.e. $\del_\a\b + \del_\b
\a = \{\a,\b\}$, for all $\a,\b \in \oneforms{M}$. Suppose $\del$
is partial, then for $\b \in \G(Q^0)$ the previous relation
becomes:
\[\del_\a \b = \{\a,\b\} = \lie{g(\a)}\b = \delta^B_{g(\a)}\b,\]
i.e. $\del$ is $Q$-adapted. Conversely, suppose $\del$ is normal
and $Q$-adapted, then $\del_\a \b = \{\a,\b\} - \del_\b \a$. Let
$\a \in \G(Q^0)$, then the right hand side of this equation is
zero, and thus $\del_\a \b =0$ for all $\a \in \G(Q^0)$ and $\b
\in \oneforms{M}$. This proves the proposition. \qed
\end{pf} We will now describe a general method for constructing normal
$g$-connections.

Let $[ \cdot\, ,\, \cdot ] : \oneforms{M} \times \oneforms{M} \to
\oneforms{M}$ denote a skew-symmetric bracket that is $\R$-linear
in both arguments and satisfies, for any $f \in {\mathcal F}(M)$,
$[\a, f \b] = g(\a)(f)\b + f[\a,\b]$. Given such a bracket on
$\oneforms{M}$, one can define a unique normal $g$-connection
$\del$ for which $[\a,\b] = \del_\a \b -\del_\b \a$, namely:
\[\del_\a \b = \half\left( [\a,\b] + \{\a,\b\}\right).\] Conversely,
given a normal $g$-connection $\del$, one can define a
skew-symmetric bracket with the desired properties by putting $
[\a,\b] = \del_\a \b -\del_\b \a $. Henceforth, we shall denote
the bracket associated with a normal $g$-connection $\del$ by $[
\a ,\b ]_\del$.

As can be easily verified, for a $g$-connection $\del$ which is
both normal and $Q$-adapted, the skew-symmetric bracket satisfies:
$[\a,\eta]_\del = \delta^B_{g(\a)} \eta$, for all $\eta \in
\G(Q^0)$ and $\a \in \oneforms{M}$. Therefore, if a Riemannian
metric $G$ is chosen, with projections $\t$ and $\t^\bot$ on
$\flat_G(Q)$ and $Q^0$ respectively, and which restricts to $h$,
this bracket takes the form:
\[[\a,\b]_\del= [\t(\a), \t(\b)]_\del + \delta^B_{g(\a)}\t^\bot(\b) -
\delta^B_{g(\b)}\t^\bot(\a).\] We only have to know the value of
the bracket acting on sections of $\flat_G(Q) \cong Q$. For
example, for the $g$-connection given by $\del_\a \b =
\del^G_{g(\a)} \t(\b) + \delta^B_{g(\a)} \t^\bot(\b)$, the
associated bracket becomes:
\[[\a,\b]_\del = \flat_G\left([g(\a),g(\b)]\right) +
\delta^B_{g(\a)}\t^\bot(\b) - \delta^B_{g(\b)}\t^\bot(\a),\] where
$[g(\a),g(\b)] = \lie{g(\a)}g(\b)$ is the usual Lie bracket on
vector fields. Note, however, that there does not seem to exist a
`natural' skew-symmetric bracket on $\oneforms{M}$, independent of
the chosen Riemannian extension $G$ of $h$, which could be used to
identify a `standard' $g$-connection which is both normal and
$Q$-adapted. One might think of imposing a metric condition in
order to completely determine such a $\del$, but the following
result tells us that it is impossible to find a $Q$-adapted
$g$-connection which is also metric.
\begin{prop}
A $Q$-adapted $g$-connection is not metric. \end{prop} \begin{pf}
Let $\del$ be $Q$-adapted $g$-connection. Suppose that $\del$
leaves $\ovl g$ invariant. This can be equivalently rewritten as
$g(\del_\a \b) = \del_\a (g(\b))$ for all $\a,\eta \in
\oneforms{M}$. Let $\eta \in \G(Q^0)$, then, since $\del$ is
$Q$-adapted this equation becomes $g(\delta^B_{g(\a)}\eta)=0$ for
all $\a \in \oneforms{M}$ and $\eta \in\G(Q^0)$. However, this is
equivalent to saying that $Q$ is involutive. Indeed, from
$g(\delta^B_{g(\a)}\eta)=0$ we have
\[0=\langle \b ,g(\delta^B_{g(\a)} \eta)\rangle  = \langle
\delta^B_{g(\a)} \eta , g(\b) \rangle = -\langle \eta,
[g(\a),g(\b)]\rangle,\] for arbitrary $\a,\b \in \oneforms{M}$ and
$\eta \in \G(Q^0)$, hence $[g(\a),g(\b)] \in \G(Q)$. \qed
\end{pf}
\section{Abnormal extremals}
\label{sec:AbnormalExtremals} For the remainder of this paper we
will always restrict ourselves to curves $c$ that are immersions,
i.e. $\dot{c}(t) \neq 0$ for all $t \in \hbox{Dom} (c)$. Such
curves can always, at least locally, be seen as (part of) an
integral curve of a smooth vector field (see e.g.\ \cite[p
28]{helgason}). Bearing this in mind, we will establish in the
present section a geometrical characterization of abnormal
extremals on a manifold with a regular, non-integrable
distribution $Q$. First, we will restrict ourselves to curves that
are integral curves of a vector field. Next, we will extend the
analysis to general continuous piecewise curves tangent to $Q$,
whose smooth parts are immersions such that they can be regarded
as a concatenation of integral curves of vector fields belonging
to $\G(Q)$.

Consider a manifold $M$ equipped with a regular distribution $Q$.
Choose an arbitrary sub-Riemannian metric $h$ (e.g.\ the
restriction of some Riemannian metric on $M$) and let $\del$ be a
fixed $Q$-adapted $g$-connection associated to the sub-Riemannian
structure $(M,Q,h)$. (From the previous section we know that such
a $g$-connection can always be found.) Suppose that $c: I
\rightarrow M$ is a curve tangent to $Q$, which is (part of) an
integral curve of a vector field $X \in \G(Q)$, defined on a
neighborhood of $\im (c)$. In particular, we have that $\dot c(t)
= X(c(t))$ for all $t \in I$. Then we know that $c$ is an abnormal
extremal if there exists a section $\eta$ of $Q^0$ along $c$ such
that $\del_\a \eta(t) = 0$ for all $t\in I$, with $\a$ a
$g$-admissible curve with base curve $c$. Let $\{\phi_s\}$ denote
the (local) flow of $X$ such that for any fixed $t \in [a,b]$,
$\phi_s(c(t)) = c(t+s)$ for all $s$ for which the right-hand side
is defined. We denote the dual of the tangent map $T\phi_s$ of
$\phi_s$ by $T^*\phi_s$, i.e. for $\a \in T^*_{\phi_s(x)}M$,
$T^*\phi_s(\a)$ is the co-vector at $x$ defined by
$T^*\phi_s(\a)(Y_x)= \a(T\phi_s(Y_x))$, for all $Y_x \in T_xM$
(with $x \in \hbox{Dom}(\phi_s)$) . We can now prove the following
lemma.
\begin{lem} Let $c: I \rightarrow M$ be an integral curve of
$X \in \G(Q)$ and let $\eta$ be an arbitrary section of $Q^0$
along $c$. Then, for any $g$-admissible curve $\a$ with base curve
$c$, the following equation holds:
\[\del_\a \eta(t) = \left.\frac{d}{ds}\right|_{s=0} \left(
T^*\phi_s (\eta(t+s))\right), \ \forall t \in I.\]
\end{lem}
\begin{pf}
Fix an arbitrary $t \in I$ and choose a local coordinate
neighborhood of $M$ containing the point $c(t)$. Since $\del_\a
\eta(t)$ is independent of the $g$-admissible curve $\a$
projecting onto $c$, we can choose $\a(t)= \ovl \a(c(t))$, where
$\ovl \a =\flat_G(X)$ and $G$ is any Riemannian metric restricting
to $h$. From Subsection
\ref{sbs:subRiemannianStructuresPreliminaryDefinitions} we know
that $g(\ovl \a)=X$, which implies indeed that $\a(t)=\ovl
\a(c(t))$ is a $g$-admissible curve with base curve $c$. In
coordinates, $\del_\a \eta(t)$ reads:
\[
\del_\a \eta(t) = \left(\dot \eta_i(t) + \fpd{g^{jk}}{x^i}(c(t))
\a_k(t)\eta_j(t)\right) dx^i|_{c(t)}.
\]
Since, for fixed $t$ and for sufficiently small $s$, the mapping
$s \mapsto T^*\phi_s (\eta(t+s))$ defines a curve in the fibre
$T_{c(t)}^*M$, the derivative at $s=0$ is well defined and can be
identified with an element of $T_{c(t)}^*M$. In coordinates this
curve is given by
\[
T^*\phi_s (\eta(t+s)) = \fpd{\phi^j_s}{x^i}(c(t))\eta_j(t+s)
dx^i|_{c(t)},
\]
and its derivative at $s=0$
equals\[\left.\frac{d}{ds}\right|_{s=0} \left( T^*\phi_s
(\eta(t+s))\right) = \left(\dot \eta_i(t) + \fpd{X^j}{x^i}(c(t))
\eta_j(t)\right) dx^i|_{c(t)}.\] Using the fact that $g(\ovl \a)=
X$ this leads to the desired result, since
\[\fpd{X^j}{x^i}(c(t))\eta_j(t) = \fpd{(g^{jk}
\ovl \a_k)}{x^i}(c(t))\eta_j(t) =
\fpd{g^{jk}}{x^i}(c(t))\a_j(c(t))\eta_j(t)\;,\] where the second
equality follows from $\eta \in \G(Q^0)$. \qed
\end{pf}
Herewith, we derive the following characterization of an abnormal
extremal.
\begin{prop} \label{propAB}
Let $c:I=[a,b] \to M$ be a curve tangent to $Q$, such that it is
an integral curve of a vector field $X \in \G(Q)$ with flow
$\{\phi_s\}$. Then, $c$ is an abnormal extremal if and only if
there exists a section $\eta$ of $Q^0$, defined along $c$, such
that $\eta(t)= T^*\phi_{-(t-a)} (\eta(a))$ for all $t \in I$.
\end{prop}
\begin{pf}
According to Theorem \ref{thm2}, $c$ is an abnormal extremal iff
there exists a section of $Q^0$ along $c$ such that $\del_\a
\eta(t) =0$, with $\a$ a $g$-admissible curve. Using the preceding
lemma, this is still equivalent to
\[
\left.\frac{d}{ds}\right|_{s=0} \left( T^*\phi_s
(\eta(t+s))\right) = 0, \ \forall t \in I.\] Acting with the map
$T^*\phi_{(t-a)} : T^*_{c(t)}M \rightarrow T^*_{c(a)}M$ on both
sides of this equation, we obtain the equivalent condition:
\[
\left.\frac{d}{dt}\right|_{t} \left( T^*\phi_{(t-a)}
(\eta(t))\right) = 0, \ \forall t \in I,\]from which it follows
that $T^*\phi_{(t-a)}(\eta(t)) = \eta(a)$. \qed
\end{pf}
This characterization of abnormal extremals that are integral
curves of a vector field, leads us to the following construction.
Let $c$ be a curve tangent to $Q$, with domain $I=[a,b]$, such
that it it is an integral curve of a vector field $X \in \G(Q)$.
For each $t \in I$ consider the subset $c^*_tQ$ of the tangent
space $T_{c(t)}M$, given by
\[c^*_tQ = \mbox{Span}\{T\phi_{-s}(Y_{c(t+s)})\ | \ \forall Y \in
Q_{c(t+s)}, \ s \in [a-t, b-t]\}.\] It is immediately verified
that $c^*_tQ$ is in fact a linear subspace of $T_{c(t)}M$.
Moreover, it is also easily seen that, for each $t \in I$ and $s
\in [a-t,b-t]$: $T\phi_s (c^*_tQ) = c^*_{t+s} Q$. Therefore, the
dimension of the linear space $c^*_tQ$ is independent of $t$. As
an aside of the following theorem it will follow that $c^*_tQ$
only depends on the set $\{\dot c(s) = X(c(s))\;|\;s \in [a,b]\}$.
\begin{thm} \label{thmAB}
Let $c$ be a curve tangent to $Q$ with domain $I=[a,b]$, such that
it is an integral curve of $X \in \G(Q)$. Then $c$ is an abnormal
extremal if and only if $c^*_aQ \neq T_{c(a)}M$.
\end{thm} \begin{pf} Suppose that $c^*_aQ \neq T_{c(a)}M$, i.e.
there exists a non zero $\eta_a \in (c^*_aQ)^0 \subset
T^*_{c(a)}M$. Define a curve $\eta$ in $T^*M$ along $c$ by
$\eta(t)= T^*\phi_{-(t-a)}(\eta_a)$. Note that $\eta(t) \neq 0$
for all $t$. We now prove that $\eta(t) \in Q^0$ and, hence,
$c(t)$ is an abnormal extremal (see Proposition \ref{propAB}). For
any $Y_{c(t)} \in Q_{c(t)}$, we have to show that $\langle
\eta(t), Y_{c(t)} \rangle =0$. By definition of $\eta(t)$ this is
indeed the case, since
\[ \langle \eta(t), Y_{c(t)} \rangle = \langle \eta_a,
T\phi_{-(t-a)}(Y_{c(t)})\rangle \mbox{ and }
T\phi_{-(t-a)}(Y_{c(t)}) \in c^*_aQ.\] Conversely, suppose that
$c(t)$ is an abnormal extremal, then, again in view of Proposition
\ref{propAB}, there exists a section $\eta$ of $Q^0$ along $c$,
which does not intersect the zero section, such that $\eta(t) =
T^*\phi_{-(t-a)}(\eta(a))$. Since $\eta(t)\in Q^0$, we then have
that $\langle \eta(t), Y_{c(t)}\rangle =0$ for all $t$ and for
arbitrary $Y_{c(t)} \in Q_{c(t)}$. This relation can be rewritten
as follows:
\[\langle \eta(t), Y_{c(t)}\rangle = \langle \eta(a),
T\phi_{-(t-a)}(Y_{c(t)}) \rangle = 0\;,\] and, hence, we conclude
that $0 \neq \eta(a) \in (c^*_aQ)^0$, which completes the
proof.\kern-.9mm \qed \end{pf}

>From the above proof it follows that each element of $(c^*_aQ)^0$
determines a unique section of $Q^0$ along $c$ by parallel
transport with respect to a $Q$-adapted $g$-connection, and vice
versa. Since for a parallel section $\eta$ of $Q^0$ along $c$ the
equation $\del_\a \eta(t) =0$ only depends on the tangent vector
to the base curve $c$, one may indeed conclude that the space
$c^*_aQ$ only depends on $\{\dot c(t)\;|; t \in [a,b]\}$.
\begin{rem}
Given a vector field $X \in \G(Q)$ such that $X(c(t)) =\dot c(t)$,
consider the subspace of $T_{c(a)}M$ spanned by $Q_{c(a)}$ and by
all tangent vectors of the form $[X,[X, \ldots [X,Y] \ldots
](c(a))$ for arbitrary $Y \in \G(Q)$, and let us denote this space
by $D_{c(a)}$. It is not difficult to prove that the space spanned
by $D_{c(a)}$ is contained in (but, in general differs from)
$c^*_aQ$.
\end{rem}
A well known result concerning abnormal extremals (see, for
instance, \cite{stri}) states that if $Q$ is `strongly bracket
generating', i.e. if $T_xM = \mbox{ Span}\{ Y(x) + [X,Y'](x) \ | \
Y,Y' \in \G(Q)\}$ for every $x \in M$ and $X \in \G(Q)$, then
there are no abnormal extremals. Since $\mbox{ Span}\{ Y(x) +
[X,Y'](x) \ | \ Y,Y' \in \G(Q)\} \subset D_{c(a)}$, the previous
remark shows that this result is compatible with Theorem
\ref{thmAB}. At least for the class of curves we are considering
here, this result can even be generalized in the following sense.
If for some $X \in \G(Q)$ we have that at every point $x \in
\hbox{Dom}(X)$ we have $T_xM = D_x$, then no integral curve of $X$
passing trough the point $x$ can be an abnormal extremal.

So far, we have only characterized those abnormal extremals that
can be regarded as integral curves of a vector field tangent to
$Q$. We shall now extend Theorem \ref{thmAB} to the class of
abnormal extremals that may be continuous piecewise curves.

Given any curve $c: I=[a,b] \to M$ tangent to $Q$ which is an
immersion, then there exists a finite subdivision of $I$, such
that the restriction of $c$ to each subinterval is an integral
curve of a vector field tangent to $Q$ (cf.\ \cite[p
28]{helgason}). This further implies that, given any continuous
piecewise curve $c : I=[a,b] \to M$ tangent to $Q$, we can apply
this property to each smooth part $c^i: [a_i,a_{i+1}]\rightarrow
M$ of $c$ (for $i = 1, \ldots, k$), where we are using the
conventions of Subsection
\ref{sbs:NormalAndAbnormalExtremalsLocalEquations}. More
precisely, each sub-curve $c^i$ can be regarded by itself as a
concatenation of integral curves of (local) vector fields
belonging to $\G(Q)$. For the sake of clarity, we will now
consider the simple case of a continuous piecewise curve
consisting of a concatenation of two integral curves of vector
fields tangent to $Q$. This will suffice to show how to proceed in
the general case of continuous piecewise curves.

Let $c: [a,b] \rightarrow M$ be a continuous piecewise curve
consisting of two smooth sub-curves $c^1: [a_1,a_2] \rightarrow M$
and $c^2: [a_2,a_3] \rightarrow M$, where $a_1=a < a_2 <a_3=b$ and
$c^i(t)=c(t)$ for $t\in ]a_i,a_{i+1}]$, and whereby we assume that
both $c^1$ and $c^2$ are integral curves of vector fields $X^1 \in
\G(Q)$ and $X^2 \in \G(Q)$, respectively. Denote the local flow of
$X^i$ by $\{\phi^i_s\}$, $i=1,2$. Since $\dot c^i (t) =
X^i(c^i(t))$ we have: $c^i(t) = \phi^i_{(t-a_i)}(c^i(a_i))$,
$i=1,2$. Consider the subspace $c^*_aQ$ of $T_{c(a)}M$ given by
\[c^*_aQ = (c^1)^*_aQ + T\phi^1_{-(a_2-a_1)}\left((c^2)^*_{a_2}
Q\right),\] where the spaces $(c^i)^*_{a_i}Q$ are defined as
above. Assume that $\eta_a \in (c^*_aQ)^0$. Then the continuous
piecewise curve in $Q^0$, defined by
\[\eta(t) =
\left\{\begin{array}{ll} T^*\phi^1_{-(t-a_1)}(\eta_a) & \ \forall t \in [a_1,a_2],\\
T^*\phi^2_{-(t-a_2)}(T^*\phi^1_{-(a_2-a_1)}(\eta_a)) & \ \forall t
\in [a_2,a_3],
\end{array}\right.\] is a parallel transported section of $Q^0$ with respect
to a $Q$-adapted connection (apply Proposition \ref{propAB} to
$c^1$ and $c^2$). This proves that if $c^*_a Q \neq T_{c(a)}M$
then $c$ is an abnormal extremal. Conversely, assume that $c$ is
an abnormal extremal. By definition there exist parallel
transported sections $\eta^1$ and $\eta^2$ of $Q^0$ along $c^1$
and $c^2$ respectively, such that $\eta^1(a_2) = \eta^2(a_2)$.
Theorem \ref{thmAB} implies that $0\neq \eta^1(a_1) \in
((c^1)^*_{a_1} Q)^0$ and $0 \neq \eta^2(a_2)\in
((c^2)^*_{a_2}Q)^0$. Since $\eta^2(a_2) = \eta^1(a_2) =
T^*\phi^1_{-(a_2-a_1)}(\eta^1(a_1))$, we conclude that
\[0 \neq \eta^1(a_1) \in \left((c^1)^*_{a_1}Q\right)^0\cap \left(
T\phi^1_{-(a_2-a_1)}((c^2)^*_{a_2} Q)\right)^0= (c^*_a Q)^0.\]
This reasoning can now be easily extended to the case where $c$ is
a general continuous piecewise curve tangent to $Q$ (for which
$\dot{c}(t) \neq 0$ at all points where the derivative is
defined). Summarizing, we have derived the following
characterization of abnormal extremals within the class of
continuous piecewise curves.

\begin{thm} Let $c:I=[a,b] \to M$ be a continuous piecewise curve
tangent to $Q$, with $\dot{c}(t) \neq 0$ at each point where the
derivative exists. Then, there always exists a finite subdivision
of $I$, with endpoints $a_1=a < a_2 < \ldots < a_{\ell} < a_{\ell
+1} =b$, such that $c$ is a concatenation of integral curves
$c^i:[a_i,a_{i+1}] \to M$ of vector fields $X^i$ tangent to $Q$,
with flow $\{\phi^i_s\}$, $i=1, \ldots , \ell$. We then have that
$c$ is an abnormal extremal if and only if
\[T_{c(a)} M \neq c^*_aQ := (c^1)^*_{a_1}Q + \sum_{i=2}^{\ell}T\phi^1_{-(a_2-a_1)}
\ldots T\phi^{i-1}_{-(a_i -a_{i-1})}((c^i)^*_{a_i}Q).\]\end{thm}
Note that, although we have used the theory of $g$-connections
associated to a sub-Riemannian structure for its derivation, the
above characterization of abnormal extremals is independent of the
choice of a sub-Riemannian metric, but only depends on the
geometry of the given distribution $Q$. This is indeed in full
agreement with the notion of abnormal extremal.
\begin{rem} While finalizing this paper,
we have come across a recent paper by P. Piccione and D.V. Tausk
\cite{piccione}, in which, following a different approach, a
similar characterization for abnormal extremals was obtained.
\end{rem} We shall now give two examples to illustrate the previous results.
\begin{exmp}
\label{sbs:Examples} {\rm Here we consider an example of abnormal
extremals, con\-struct\-ed by R. Montgomery \cite{Mont}. Let $M=
\R^3-\{0\}$ and let $Q$ be the 2-dimensional distribution spanned
by the vector fields (expressed in cylindrical coordinates): $X_1
= \fpd{}{r}, X_2 = \fpd{}{\theta} - F(r) \fpd{}{z}$, where $F(r)$
is a function on $M$ with a single non degenerate maximum at
$r=1$, i.e. $F$ satisfies: \[ \left. \frac{d}{dr} F(r)
\right|_{r=1} = 0 \quad \hbox{\rm and} \quad \left.
\frac{d^2}{dr^2} F(r) \right|_{r=1} < 0.\]Such a function can
always be constructed (take, for instance, $F(r) = \half r^2 -
\frac{1}{4}r^4$). The distribution thus defined is everywhere of
rank two, and is differentiable by definition. The flows of
$X_1,X_2$ are denoted by $\{\phi_s\}$, $\{\psi_s\}$, respectively.
In particular, we have $\phi_t(r,\theta,z)= (t+r,\theta,z)$,
$\psi_t(r,\theta,z)= (r,\theta +t,z-F(r)t)$. Let $c:[0,1]
\rightarrow M$ be an integral curve of $X_1$ through
$x_0=(r_0,\theta_0,z_0)$ at $t=0$. The subspace
\[ c^*_0Q = \mbox{Span}\left\{X_1(x_0), X_2(x_0),
\left. \fpd{}{\theta}\right|_{x_0}
-F(r+t)\left.\fpd{}{z}\right|_{x_0} \left.\right| \forall t \in
[0,1]\right\}.\]This subspace coincides with the whole tangent
space at $x$, as can be seen from: \bea
v_r\left.\fpd{}{r}\right|_{x_0} + v_\theta
\left.\fpd{}{\theta}\right|_{x_0} + v_z
\left.\fpd{}{z}\right|_{x_0} &=& v_rX_1(x_0) + v_\theta X_2(x_0)\\
&& \mbox{} + \frac{v_z + v_\theta F(r_0)}{F(r_0+t) -
F(r_0)}(X_2-\phi^*_tX_2)(x_0),\eea where $t$ is chosen such that
$F(r_0+t)\neq F(r_0)$. So, in view of Theorem \ref{thmAB}, one can
conclude that an integral curve of $X_1$ can not be an abnormal
extremal. Let $c': [0,1] \rightarrow M$ be an integral curve of
$X_2$, with $c'(0) = x_0=(r_0,\theta_0,z_0)$. Then we have
\[{c'}_0^*Q = \mbox{Span}\left\{ X_1(x_0), X_2(x_0),
\left.\fpd{}{r}\right|_{x_0} +
F'(r_0)t\left.\fpd{}{z}\right|_{x_0} \ | \ \forall t \in
[0,1]\right\}.\] If $x_0$ is a point on the cylinder defined by
$r=1$, then one easily sees that $c'^*_0Q \neq T_xM$ since
$F'(1)=0$. Therefore, every helix $c': [0,1] \to \R^3: t \mapsto
(1, \theta + t , z-F(1)t)$ is an abnormal extremal, i.e. there
exists a section of $Q^0$ along the curve $c'$ through
$x_0=(1,0,0)$ such that
\[\eta(t):=T^*\psi_{-t}(F(1)\left.d\theta\right|_x +
\left.dz\right|_x) = F(1)\left.d\theta\right|_{(1,t,-F(1)t)} +
\left.dz\right|_{(1,t,-F(1)t)}.\]} \end{exmp}
\begin{exmp}{\rm
We now treat an example that was constructed by W. Liu and H.J.
Sussmann, \cite{sussmann2}. Let $M=\R^3$ and $Q$ spanned by $X_1=
\fpd{}{x}, X_2= (1-x)\fpd{}{y} + x^2\fpd{}{z}$, where we use
cartesian coordinates, $x,y,z$. The flows $\{\phi_s\}$ of $X_1$
and $\{\psi_s\}$ of $X_2$ are given by $\phi_t(x,y,z)=(x+t,y,z)$
and $\psi_t(x,y,z)= (x, (1-x)t +y, x^2t +z)$.  The pull-back of
$X_1$ under $\psi_t$ equals $\psi^*_tX_1=\fpd{}{x} + t\fpd{}{y} -
2xt \fpd{}{z}$, and this vector field can be written as a linear
combination of $X_1,X_2$ for any value of $t$ and at all points
for which $x=0$ or $x=2$. Indeed, if $x=0$, then
$\psi^*_tX_1(0,y,z) = X_1(0,y,z) + tX_2(0,y,z)$. If $x=2$, then
$\psi^*_tX_1(2,y,z) = X_1(2,y,z) -t X_2(2,y,z)$. Therefore, each
curve defined by $c: I \rightarrow M: t\mapsto (x, (1-x)t +y,
x^2t+z)$ for any given point $(x,y,z)$ with $x=0$ or $x=2$, is an
abnormal extremal.}
\end{exmp}
To end this section, we present a construction for the tangent
vector to certain variations of a given curve $c: [a,b] \to M$
tangent to $Q$, that have been used in a derivation of the Maximum
principle in \cite{Pont}. We shall see that the set of all such
tangent vectors determines a subspace of the tangent space $T_bM$
that equals $c^*_bQ$. Suppose that $c: [a,b] \rightarrow M$ is a
curve tangent to $Q$, which is an integral curve of a vector field
with flow $\{\phi_t\}$, such that $c(a+t)= \phi_t(c(a))$. The type
of variations of $c$ we have in mind here, are specified by a
triple $(Y, \t, \delta t)$ with $Y \in \G(Q)$, $\t \in [a,b]$ and
$\delta t \ge 0 \in \R$. Denote the flow of $Y$ by $\{\psi_s\}$.
The variation $\tilde c : [a,b] \times \R \to M$, associated to
the triple $(Y,\t,\delta t)$ for $\t \in]a,b]$, is then defined
by:\[\tilde c(t,\epsilon) = \left\{\begin{array}{ll} c(t) & \ \ a
\le t \le \t-\epsilon \delta t,\\ \psi_{t-(\t -\epsilon \delta
t)}(c(\t-
\epsilon \delta t)) &\ \ \t - \epsilon \delta t \le t \le \t, \\
\phi_{t - \t}(\psi_{\epsilon \delta t}(c(\t -\epsilon \delta t)))
&\ \ \t \le t \le b,
\end{array}\right. \]which is well defined for $\epsilon$ small
enough. For $\t =a$, a slightly different definition for $\tilde c
: [a, b] \to M$ is needed: $\tilde c(t,\epsilon )=
\phi_{t-a}(\psi_{\epsilon \delta t}(\phi_{-\epsilon \delta
t}(c(a))))$. The tangent vector to any variation $\tilde c$ at
$(t, \epsilon) = (b,0)$ equals:
\[V(Y,\t, \delta t) = T\phi_{b-\t}(\delta t Y(c(\t)) - \delta t
\dot c(\t)).\] Since $Y(c(\t)) - \dot c(\t) \in Q_{c(\t)}$, the
vector $V(Y,\t, \delta t)$ belongs to $c^*_bQ$. Even more, the
space spanned by all $V(Y,\t, \delta t)$ with $Y \in \G(Q)$, $\t
\in [a,b]$ and $\delta t \in \R$, equals $c^*_bQ$. Therefore, the
necessary and sufficient condition from Theorem \ref{thmAB}
measures the dimensionality of the space spanned by tangent
vectors to variations. A more detailed discussion will be
presented in a forthcoming paper in which we will construct a
natural connection over a bundle map associated with a control
problem, which will lead to a weaker version of the Maximum
principle.
\section{Normal extremals}
\label{sec:NormalExtremals} In this section we will make use of
Theorem \ref{thm1} to recover some known results about normal
extremals. Consider a sub-Riemannian structure $(M,Q,h)$ and let
$G$ be an arbitrary Riemannian metric on $M$ restricting to $h$.
Theorem \ref{thm1} then says $\del_{\a}\a(t) = \del^G_{\dot
c}\t(\a)(t) + \delta^B_{\dot c}\t^\bot(\a)(t),$ where $\a$ is a
$g$-admissible curve with base curve $c$, and $\del$ is any normal
$g$-connection. This immediately leads to the following result.
\begin{prop}
Let $c:I \rightarrow M$ be a curve tangent to $Q$ that is a
geodesic with respect to a Riemannian metric $G$ restricting to
$h$, then $c$ is a normal extremal.
\end{prop}
\begin{pf}
The curve $c$ is a normal extremal if there exists a
$g$-admissible curve $\a$ with base curve $c$, which is
auto-parallel with respect to a normal $g$-connection $\del$.
Since $c:I \rightarrow M$ is a geodesic with respect to $G$, i.e.
$\del^G_{\dot c} \dot c(t)=0$ $\forall t \in I$, we know from
Section \ref{sbs:subRiemannianStructuresPreliminaryDefinitions}
that $\a=\flat_G(c)$ is a $g$-admissible curve with base curve $c$
for which $\t(\a) = \a$ or  $\t^\bot(\a)=0$. It then follows that
$\del_\a \a(t) =0$ since $\del_\a \a(t) = \del^G_{\dot c}
\t(\a)(t) = \flat_G(\del^G_{\dot c} \dot c(t)) =$0.\qed\end{pf}
Let $c: I =[a,b] \rightarrow M$ be a normal extremal. Then there
exists a $g$-admissible curve $\a$ which is auto-parallel with
respect to a normal $g$-connection. Given any $t_0 \in I$, then
one can always find a one form $\ovl \a$ and a compact subinterval
$J$ of $I$ containing $t_0$, such that $\ovl \a (c(t)) = \a(t)$
for all $t \in J$ and $c(J)$ is contained in a coordinate
neighborhood $U$. We will now construct a local Riemannian metric
$G$ restricting to $h$ on $Q$ such that $c|_J$ is a geodesic with
respect to this Riemannian metric.

Since $g(\ovl \a) \neq 0$, one can construct a local basis of
$\oneforms{U}$, namely $\{ \ovl \a = \b^1, \b^2, \ldots , \b^n\}$,
such that $\b^{k+1}, \ldots , \b^{n}$ determine a local basis for
$\G(Q^0)$, defined on $U$. Let $\{X_1, \ldots , X_n\}$ denote the
dual basis of $\vectorfields{U}$. Then the vector fields $X_j$,
for $j=1, \ldots ,k$, form a local basis for $\G(Q)$, since
$\langle \b^i,X_j\rangle \equiv 0$ for $i=k+1, \ldots, n$. We can
now define a Riemannian metric $G$ on $U$, restricting to $h$, as
in Section
\ref{sbs:subRiemannianStructuresPreliminaryDefinitions}, i.e.\ for
arbitrary vector fields $Y$ and $Z$ on $U$,
\[G(x)(Y,Z) = \sum_{r,s=1}^k Y^rZ^sh(x)(X_r(x),X_s(x)) + \sum_{r=k+1}^n Y^rZ^r,\]
where we have put $Y(x) = Y^rX_r(x)$ and $Z(x) = Z^rX_r(x)$ for
some $Y^r,Z^r \in \R$ ($r=1,\ldots ,n)$. From the definition of
$G$ we can derive that $Q^\bot$ is spanned by $\{X_{k+1}, \ldots ,
X_n\}$ or $\t^\bot(\ovl \a) = 0$, implying that $\t^\bot(\a(t))=0$
or $\flat_G(\dot c(t)) = \a(t)$. From $\del_\a \a (t) = 0$ and
$\t^\bot(\a(t))=0$ we obtain $\del^G_{\dot c} \dot c(t) =0$ for
any $t \in J$.

\begin{prop}
Let $c:I \to M$ be a normal extremal. Then for any $t \in I$ there
exists a compact neighborhood $J$ of $t$ such that $c$ restricted
to $J$ is a geodesic with respect to some Riemannian metric
restricting to $h$ on $Q$.
\end{prop}
This proves, in particular, that a normal extremal is locally
length minimizing.

Let $c$ be a normal extremal and let $\del$ be a normal and
$Q$-adapted $g$-connection (recall that such a $\del$ always
exists). Suppose that $c$ is degenerate in the following sense:
there exist two $g$-admissible curves $\a,\b$ with base curve $c$,
such that $\del_\a \a(t)= \del_\b \b(t) =0$. We will now see that
$c$ is then also an abnormal extremal. We have proven before that
a normal and $Q$-adapted connection is partial, i.e. $\del_\a
=\del_\b$ if $g(\a)=g(\b)$. Therefore one obtains that $\del_\a(\a
-\b)(t)=0$. Since $g(\a(t)-\b(t))=0$, or $\eta(t)= (\a-\b)(t) \in
Q^0$ for all $t$, $\eta$ is a parallel transported section along
$\a$, lying entirely in $Q^0$ and, hence, $c$ is an abnormal
extremal. Conversely, assume that $c$ is a normal extremal, i.e.
$c$ is the base curve of an auto-parallel curve $\a$ with respect
to $\del$, and that $c$ is also an abnormal extremal. Let $\eta$
denote a parallel transported section along $\a$ lying in $Q^0$.
Then, using the same arguments as before, $\a + \eta$ is also an
auto-parallel curve with base curve $c$. We can conclude that
curves that are both normal and abnormal are degenerate in the
sense that they admit more than one $g$-admissible curve that is
auto-parallel.

\section{Vakonomic dynamics and nonholonomic mechanics}
\label{sec:VakonomicDynamicsAndNonholonomicMechanics} As a natural
consequence of the approach to sub-Riemannian structures in terms
of generalized connections, we will see how to establish
coordinate independent conditions for the motions of a free
mechanical system subjected to linear nonholonomic constraints to
be normal extremals with respect to the associated sub-Riemannian
structure, and vice versa. We first give a definition of what we
understand under a free mechanical systems subjected to linear
nonholonomic constraints (shortly free nonholonomic mechanical
system) and the associated sub-Riemannian structure.

Assume that a manifold $M$ is equipped with a non-integrable
regular distribution $Q$ on $M$ and a Riemannian metric $G$. A
free mechanical system with linear nonholonomic constraint $Q$
consists of a free particle with Lagrangian $L(v)= \half G(v,v)
\in {\cal F}(TM)$, subjected to the constraint $v \in Q$. (``Free"
refers here to the absence of external forces.) The problem of
determining the \emph{dynamics of the free nonholonomic mechanical
system} then consists in finding the solutions of the following
equation (see \cite{bloch,frans})\[\pi(\del^G_{\dot c}\dot c(t))=0
\quad \mbox{ and } \quad \dot c(t) \in Q \;, \forall t,\]where
$\pi$ is the orthogonal projection of $TM$ onto $Q$ with respect
to $G$ and $\del^G$ the Levi-Civita connection associated with
$G$. The associated sub-Riemannian structure is given by
$(M,Q,h_G)$, with $h_G$ the restriction of $G$ to $Q$.

In \cite{mijzelf} we have constructed a unique generalized
connection $\del^{nh}$ over the bundle map $i: Q \hookrightarrow
TM$ on the linear bundle $Q$, namely: $\del^{nh}_XY = \pi(
\del^G_X Y)$ (we have identified $X \in \G(Q)$ with $i \circ X \in
\vectorfields{M}$). The $i$-connection $\del^{nh}$ preserves the
sub-Riemannian metric $h_G$ on $Q$, i.e. $\del^{nh}_X h_G = 0$ for
any $X \in \G(Q)$, and satisfies $\del^{nh}_XY - \del^{nh}_Y X -
\pi[X,Y] =0$ for all $X,Y \in \G(Q)$. One can prove that
$\del^{nh}$ is completely determined by these two properties. In
this setting, the $i$-admissible curves are precisely curves
tangent to $Q$. Therefore, a motion $c$ of the free nonholonomic
mechanical system is characterized by the condition that
$\del^{nh}_{\dot c}\dot c(t)= 0$, for all $t$.

The \emph{vakonomic dynamical problem}, associated with the free
particle with linear nonholonomic constraints, consists in finding
normal extremals with respect to the associated sub-Riemannian
structure $(M,Q,h_G)$. It is interesting to compare the solutions
of the nonholonomic mechanical problem with the solutions of the
vakonomic dynamical problem, because the equations of motion for
the mechanical problem are derived by means of d'Alembert's
principle, whereas the normal extremals are derived from a
variational principle. This has been discussed for more general
Lagrangian systems by J. Cort\' es, \etal{} \cite{manolo}. For the
free particle case, we shall present here an alternative
(coordinate free) approach .
\begin{defn}
Given a Riemannian metric $G$ and a regular distribution $Q$ on a
manifold $M$. We can then define the following two tensorial
operators:\[
\begin{array}{l} \Pi^G: \G(Q)\otimes \G(Q) \rightarrow \G(Q^\bot),
(X,Y) \mapsto \pi^\bot(\del^G_XY),\\ \Pi^B: \G(Q) \otimes \G(Q^0)
\rightarrow \G(({Q^\bot})^0), (X,\eta) \mapsto
\t(\delta^B_X\eta).\end{array}\]
\end{defn}
It is indeed easily seen that both $\Pi^G$ and $\Pi^B$ are ${\cal
F}(M)$-bilinear in their arguments and, hence, their action can be
defined point-wise, with expressions like $\Pi^G(X_x,Y_x)$ and
$\Pi^B(X_x,\eta_x)$, for $X_x,Y_x \in Q_x$ and $\eta_x \in Q^0$,
having an obvious and unambiguous meaning.

The operator $\Pi^B$ is related to the `curvature' of the
distribution $Q$ as follows: let $X,Y \in \G(Q)$, then one has:
\[\langle \Pi^B(X,\eta) , Y \rangle = \langle \delta^B_X\eta, Y
\rangle = - \langle \eta, [X,Y] \rangle, \mbox{ for any } \eta \in
\G(Q^0). \] Thus $\Pi^B \equiv 0$ if and only if $Q$ is
involutive. The following lemma shows the importance of these
tensors. First, define a linear connection $\widetilde{\del}^B$
over $i: Q \hookrightarrow TM$ on the bundle $Q^0$ by the
prescription $\widetilde{\del}^B_X \eta = \t^\bot (\delta^B_X
\eta)$ with $X \in \G(Q)$ and $\eta \in \G(Q^0)$.
\begin{lem}
Given a Riemannian metric $G$ and a regular distribution $Q$ on a
manifold $M$. Assume that $c: I=[a,b] \rightarrow M$ is a curve
tangent to $Q$ and let $\del$ be a $Q$-adapted $g$-connection with
respect to the associated sub-Riemannian structure $(M,Q,h_G)$.
Then, the following properties hold:
\begin{enumerate}
\item \label{tensor1}
Given $Y_a \in Q_{c(a)}$, denote the parallel transported curves
along $c$, with initial point $Y_a$, with respect to $\del^{nh}$,
resp.\ $\del^G$, by $\tilde Y(t)$, resp.\ $Y(t)$. Then $\tilde Y
(t) = Y(t)$ for all $t$, if and only if $\Pi^G(\dot c(t), \tilde
Y(t)) = 0$ for all $t \in I$.
\item \label{tensor2}
Given $\eta_a \in Q^0_{c(a)}$, denote the parallel transported
curves along $c$, with initial point $\eta_a$, with respect to
$\widetilde{\del}^B$, resp.\ $\del$, by $\tilde \eta(t)$, resp.\
$\eta(t)$.  Then $\tilde \eta(t) = \eta(t)$ if and only if
$\Pi^B(\dot c(t), \tilde \eta(t))=0$.
\end{enumerate}
\end{lem}
\begin{pf} {\it (\ref{tensor1})}
>From the definition of $\Pi^G$ it follows that, given any section
$\tilde Z(t)$ of $Q$ along $c$, the following equation holds:
$\del^{nh}_{\dot c} \tilde Z(t) = \del^G_{\dot c} \tilde Z(t) -
\Pi^G(\dot c(t), \tilde Z(t))$. Assume that $\tilde Z(t) = \tilde
Y(t) = Y(t)$, then we have $\Pi^G(\dot c(t), \tilde Y(t))=0$. This
already proves the statement in direction. The converse follows
from the fact that parallel transported curves with respect to any
connection are uniquely determined by their initial conditions.

The proof of {\it(\ref{tensor2})} follows from similar
arguments.\qed
\end{pf} Note that property {\it (\ref{tensor2})} of the previous lemma
gives necessary and sufficient conditions for the existence of
curves that are abnormal extremals, i.e.: $c$ is an abnormal
extremal if and only if there exists a parallel transported
section $\tilde \eta$ of $Q^0$ along $c$ with respect to
$\widetilde{\del}^B$ such that, in addition, $\Pi^B(\dot c(t),
\tilde \eta(t)) = 0$ for all $t$. We shall now investigate some
further properties of the operators $\Pi^B$ and $\Pi^G$.
\begin{defn}
For $x \in M$, let $X_x$ be a non-zero element of $Q_x$. Define a
subspace of $T_xM$ as follows: \bea Q_x + [X, Q_x] &=&
\mbox{Span}\{ Y(x) + [\tilde X,Y'](x) \ |\ Y,Y' \in \G(Q);\\&&
\mbox{} \tilde X \in \G(Q) \mbox{ with } \tilde X(x) =
X\}.\eea\end{defn} As a side result of the following lemma, it
will be seen that the space $Q_x+[X,Q_x]$ is independent of the
extension $\tilde X$ of $X_x$ used in its definition and, hence,
also justifies the notation.
\begin{lem}
Let $\eta_x \in Q_x^0$ and $X_x \in Q_x$ for some $x \in M$. Then
$\Pi^B(X_x, \eta_x) = 0$ if and only if $\eta \in (Q_x +
[X,Q_x])^0$.\end{lem}
\begin{pf}
Let $\Pi^B(X, \eta)=0$. This is equivalent to $\langle \eta ,
[\tilde X, Y' ](x) \rangle =0$ for any $\tilde X,Y' \in \G(Q)$
with $\tilde X(x) = X_x$. Since $\eta_x \in Q^0_x$, we may
conclude that $\eta_x \in (Q_x + [X,Q_x])^0$. The converse follows
by reversing the previous arguments. \qed
\end{pf}
Another useful property is given by the following lemma.
\begin{lem}
Let $M$ be a manifold with a Riemannian metric $G$ and a regular
non-integrable distribution $Q$, and consider the associated
sub-Riemannian structure $(M,Q,h_G)$. Let $\del$ be a normal
$g$-connection. We then have for $\a \in \oneforms{M}$ that
$\del_\a \a = 0$ if and only if
\[\begin{array}{l}
\flat_G(\del^{nh}_{g(\a)}g(\a)) = - \Pi^B(g(\a), \t^\bot(\a)) \mbox{ and}\\
\widetilde{\del}^B_{g(\a)}\t^\bot(\a) = -
\flat_G(\Pi^G(g(\a),g(\a))).\end{array}\]
\end{lem}
\begin{pf}
>From Theorem \ref{thm1} one has that $\del_\a \a =0$ if and only
if $\del^G_{g(\a)} \t(\a) + \delta^B_{g(\a)} \t^\bot(\a) =0$.
Using the following relations \[\begin{array}{l} \t(\a) =
\flat_G(g(\a)),\\ \del^G \circ \flat_G = \flat_G \circ \del^G,\\
\del^G_{g(\a)} g(\a) = \del^{nh}_{g(\a)} g(\a) +
\Pi^G(g(\a),g(\a)),\\ \delta^B_{g(\a)} \t^0(\a) =
\widetilde{\del}^B_{g(\a)} \t^0(\a) + \Pi^B(g(\a),\t^0(\a)),
\end{array}\] together with the fact that $T^*M= \flat_G(Q) \oplus
Q^0$ and $Q^0 \cong \flat_G(Q^\bot)$, the equivalence is
immediately proven.\qed
\end{pf}
The previous lemmas can now be used to derive necessary and
sufficient conditions for a motion of a free nonholonomic
mechanical system to be normal extremals and vice versa. Let $M$
again be a manifold with a Riemannian metric $G$ and a regular
non-integrable distribution $Q$.
\begin{prop}
A solution $c: [a,b] \rightarrow M$ of a free nonholonomic system
determined by the triple $(M,Q,G)$ is a solution of the
corresponding vakonomic problem, and vice versa, if and only if
there exists a section $\eta$ of $Q^0$ along $c$ such that
\begin{equation} \label{vaknonh}
\widetilde{\del}^B_{\dot c}\eta(t) = - \flat_G(\Pi^G(\dot
c(t),\dot c(t)))
\end{equation}
and such that, in addition $\eta(t) \in (Q_{c(t)} + [\dot c(t),
Q_{c(t)}])^0$ for all $t$. \end{prop}
\begin{pf}
The condition for any $g$-admissible curve $\a(t)= \flat_G(\dot
c(t)) + \eta(t)$ with base curve $c$ (where $\eta(t)$ is any
section of $Q^0$ along $c$) to be parallel transported with
respect to a normal $g$-connection is that $\del_\a \a(t)=0$. This
can equivalently be written as:
\[\begin{array}{l}
\flat_G(\del^{nh}_{\dot c}\dot c(t)) = - \Pi^B(\dot c(t), \eta(t)) \mbox{ and}\\
\widetilde{\del}^B_{\dot c}\eta(t) = - \flat_G(\Pi^G(\dot
c(t),\dot c(t))).\end{array}\] Thus $\del^{nh}_{\dot c}\dot c(t)
=0$ if and only if $\Pi^B(\dot c(t), \eta(t)) =0$, where $\eta(t)$
is a solution of $\widetilde{\del}^B_{\dot c}\eta(t) = -
\flat_G(\Pi^G(\dot c(t),\dot c(t)))$.\qed \end{pf}
\begin{rem}
Given any $\eta_0$ in $(Q_{c(a)} + [\dot c(a), Q_{c(a)}])^0$ then
(\ref{vaknonh}) always admits a solution, $\eta(t)$ with initial
condition $\eta(a) = \eta_0$. The obstruction for $c$ to be
simultaneously a motion of the nonholonomic mechanical system and
a solution to the vakonomic dynamical problem, lies in the fact
that $\eta(t)$ should belong to $(Q_{c(t)} + [\dot c(t),
Q_{c(t)}])^0$ for all $t$, and this is not guaranteed by the fact
that $\eta(t)$ is a solution of (\ref{vaknonh}). The search for
geometric conditions for solutions $\eta(t)$ of this equation to
remain in $(Q_{c(t)} + [\dot c(t), Q_{c(t)}])^0$ for all $t$, is
left for future work.
\end{rem}
\begin{ack}
This work has been supported by a grant from the ``Bijzonder
onderzoeksfonds'' of Ghent University. Special thanks goes to F.
Cantrijn for useful discussions and support.
\end{ack}


\begin{thebibliography}{99}
\bibitem{bloch}
A.M. Bloch and P.E. Crouch, Newton's law and integrability of
nonholonomic systems, {\em Siam J. Control Optim.} {\bf 36} (1998)
2020-2039.
\bibitem{Brickell} F. Brickell and R.S. Clark, {\em Differentiable
Manifolds. An Introduction} (Van Nostrand Reinhold, London, 1970).
\bibitem{Bott}
R. Bott, S. Gitler and I.M. James, {\em Lectures on Algebraic and
Differential topolgy} (Springer, Lecture Notes in Mathematics 279,
1972).
\bibitem{onszelf}
F. Cantrijn and B. Langerock, Generalised connections over a
vector bundle map, {\em Diff. Geom. Appl.}, to appear.
\bibitem{frans}
F. Cantrijn, J. Cort\'es, M. de L\'eon and M. Mart\'\i n de Diego,
On the geometry of generalized Chaplygin system, {\em Math. Proc.
Camb. Phil. Soc.} to appear.
\bibitem{chow} W.L. Chow, \"Uber Systeme Von Linearen Partiellen
Differentialgleichungen erster Ordnung, {\em Math. Ann.}, {\bf
117} (1939) 98-105.
\bibitem{manolo}
J. Cort\' es, M. de L\'eon, D. Mart\'\i n De Diego and S. Mart\'\i
nez, Geometric description of vakonomic and nonholonomic dynamics.
Comparison of solutions, submitted for publication, Preprint:
math.DG/0006183.
\bibitem{ferdi} R.L. Fernandes, Connections in
Poisson Geometry I: Holonomy and invariants, {\em J. Diff. Geom.},
 {\bf 54} (2000) 303-365.
\bibitem{ferdi2} R.L.
Fernandes, Lie Algebroids, Holonomy and Characteristic Classes,
{\em Advances in Mathematics}, {\bf 170} (2002) 119-179.
\bibitem{helgason} S. Helgason, {\em Differential geometry, Lie groups, and symmetric
spaces}, (Academic Press, New York, 1978).
\bibitem{koba} S. Kobayashi and
K. Nomizu, {\em Foundations of differential geometry Volume I and
II}  (Intersience Publishers, London, 1963).
\bibitem{mijzelf}
B. Langerock, Nonholonomic mechanics and connections over a bundle
map, {\em J. Phys. A: Math. Gen.}, {\bf 34} (2001) L609-L615.
\bibitem{Lib}
P. Libermann and C.-M. Marle, {\em Symplectic Geometry and
Analytical Mechanics} (Reidel, Dortrecht, 1987).
\bibitem{sussmann2}
W. Liu and H.J. Sussmann, Shortest paths for sub Riemannian
metrics on rank two distributions, {\em Memoirs AMS} {\bf 118}
(1995).
\bibitem{Mont}
R. Montgomery, Abnormal Minimizers, {\em Siam Journal on control
and optimization}, {\bf 32} (1994) 1605-1620.
\bibitem{piccione} P.
Piccione and D.V. Tausk, Variational aspects of the geodesics
problem in sub-Riemannian geometry, {\em J. Geom. Phys.} {\bf 39}
(2001) 183-206.
\bibitem{Pont}
L.S. Pontryagin, V.G. Boltyanskii, R.V. Gamklelidze and E.F.
Mishchenko, {\em The Mathematical Theory of Optimal Processes},
(Wiley, Interscience, 1962).
\bibitem{stri}
R.S. Strichartz, sub-Riemannian geometry, {\em J. Diff. Geom.}
{\bf 24} (1986) 221-263.
\bibitem{stri2}
R.S. Strichartz, Corrections to ``sub-Riemannian geometry'', {\em
J. Diff. Geom.} {\bf 30} (1989) 595-596.
\bibitem{sussmann3} H.J. Sussmann, A cornucopia of
four-dimensional abnormal sub-Riemannian minimizers, in: A.
Bellaïche and J.-J. Risler, eds., {\em Sub-Riemannian geometry},
(Progr. Math., 144, Birkh\"{a}user, Basel, 1996) 341-364.
\bibitem{sussmann4} H.J.
Sussmann, An introduction to the coordinate-free maximum
principle, in: B. Jakubczyk and W. Respondek, eds., {\em Geometry
of Feedback and Optimal Control}  (Marcel Dekker, New York, 1997)
463-557.
\end{thebibliography}
\end{document}